\documentclass[a4paper,twoside,11pt]{amsart}

\usepackage{fullpage}
\usepackage{amsmath} \usepackage{amssymb} \usepackage{amsopn}
\usepackage{amsthm} \usepackage{amsfonts} \usepackage{mathrsfs}

\usepackage[dvips]{graphicx} \usepackage[usenames]{color}

\usepackage[all]{xy}

\usepackage[latin1]{inputenc}  %Umlaute im Text verwenden

\usepackage[pdfpagelabels, dvips]{hyperref}

\newcommand{\Kurztitel}{Cuspidal sections} 
\hypersetup{
  pdftitle={},
  pdfauthor={},
  pdfsubject={},
  pdfkeywords={},
  colorlinks=false,
  breaklinks=true,
  bookmarksopen=true,
  bookmarksnumbered=true,
  pdfpagemode=UseOutlines,
  plainpages=false}

\DeclareMathOperator{\rD}{D}

\DeclareMathOperator{\rH}{H}
\DeclareMathOperator{\rI}{I}

\DeclareMathOperator{\rT}{T}

\DeclareMathOperator{\rc}{c}

\DeclareMathOperator{\rh}{h}

\newcommand{\bC}{{\mathbb C}}

\newcommand{\bG}{{\mathbb G}}

\newcommand{\bL}{{\mathbb L}}

\newcommand{\bN}{{\mathbb N}}

\newcommand{\bP}{{\mathbb P}}
\newcommand{\bQ}{{\mathbb Q}}
\newcommand{\bR}{{\mathbb R}}

\newcommand{\bZ}{{\mathbb Z}}

\newcommand{\cD}{{\mathscr D}}

\newcommand{\cH}{{\mathscr H}}

\newcommand{\cM}{{\mathscr M}}

\newcommand{\dO}{{\mathcal O}}

\newcommand{\surj}{\twoheadrightarrow} 
\newcommand{\inj}{\hookrightarrow}

\DeclareMathOperator{\pr}{pr}

\DeclareMathOperator{\Hom}{Hom}

\DeclareMathOperator{\Aut}{Aut}

\DeclareMathOperator{\Spec}{Spec}

\DeclareMathOperator{\Div}{Div}
\DeclareMathOperator{\divisor}{div}
\DeclareMathOperator{\Pic}{Pic}

\newcommand{\OO}{\dO}

\DeclareMathOperator{\Alb}{Alb}

\newcommand{\Gm}{\bG_m}

\DeclareMathOperator{\Br}{Br}

\DeclareMathOperator{\Ext}{Ext} 

\DeclareMathOperator{\cHom}{\cH\mathit{om}}

\DeclareMathOperator{\res}{res}

\DeclareMathOperator{\pig}{\overline{\pi}}

\DeclareMathOperator{\Gal}{Gal}

\newcommand{\ep}{\varepsilon}
\newcommand{\ph}{\varphi}

\newcommand{\alg}{{\rm alg}}
 
\newcommand{\sh}{{\rm sh}}

\newcommand{\tors}{{\rm tors}}
\newcommand{\ab}{{\rm ab}}

\newcommand{\ty}{{\tilde{y}/y}}
\DeclareMathOperator{\cp}{{\rm Cusps}}
\DeclareMathOperator{\tcp}{{\rm \widetilde{Cusps}}}
\newcommand{\solv}{{\rm solv}}
\newcommand{\cyl}{{\rm cyl}}
\newcommand{\orient}{{\rm Or}}
\DeclareMathOperator{\cl}{cl}
\DeclareMathOperator{\Sym}{Sym}
\newcommand{\loc}{{\rm loc}}
\newcommand{\kum}{{\rm Kum}}
\newcommand{\exc}{{\rm exc}}
\newcommand{\comp}{{\rm comp}}
\DeclareMathOperator{\DV}{{\rm DV}}

\newcommand{\bruch}[2]{\genfrac{}{}{0.5pt}{}{#1}{#2}}

\newtheorem{thm}{Theorem}

\newtheorem{prop}[thm]{Proposition}
\newtheorem{lem}[thm]{Lemma}

\newtheorem{cor}[thm]{Corollary}
\newtheorem{conj}[thm]{Conjecture}

\theoremstyle{definition}

\theoremstyle{remark}
\newtheorem{rmk}[thm]{Remark}

\newenvironment{pro}[1][Proof]{{\it{#1:}} }{\hfill $\square$}
\newenvironment{pro*}[1][Proof]{{\it{#1:}} }{}

\newcounter{absatzcounter}[section]
\numberwithin{equation}{section}

\begin{document}

\hrule
\vspace{1.5cm}

\title[\Kurztitel]{On cuspidal sections of  algebraic \\ fundamental groups} 
\author{Jakob Stix}
%\subjclass[2000]{Primary ; Secondary }
\thanks{The author acknowledges support provided by DFG grant 
STI576/1-(1+2).}
\address{Mathematisches  Institut, Universit\"at Bonn, Beringstra\ss e 1, 53115 Bonn, Germany}
\email{stix@math.uni-bonn.de}
%\urladdr{http://}
\date{\today} %\date{15 Juli 2008}

\maketitle

\begin{quotation} 
  \noindent \small {\bf Abstract} --- Rational points in the boundary of a hyperbolic curve over a field with sufficiently nontrivial Kummer theory are the source for an abundance of sections of the fundamental group exact sequence. We follow and refine Nakamura's approach towards these boundary sections.
  For example, we obtain a weak anabelian theorem for hyperbolic genus $0$ curves over quite general fields including for example $\bQ^\ab$.
\end{quotation}

%%%%%%%%%%%%%%%%%%%%%%%%%%%%%%%%%%%%%%%%%%%
%%%%%%%%%%%%%%%%%%%%%%%%%%%%%%%%%%%%%%%%%%%

\setcounter{tocdepth}{1} {\scriptsize \tableofcontents}

%%%%%%%%%%%%%%%%%%%%%%%%%%%%%%%%%%%%%%%%%%%
%%%%%%%%%%%%%%%%%%%%%%%%%%%%%%%%%%%%%%%%%%%

\section{Introduction}

\subsection{The section conjecture} 
The \'etale fundamental group of a geometrically connected variety $U$ over a field $k$ with fixed algebraic closure  $k^\alg$ sits naturally in a short exact sequence 
\[ 1 \to \pi_1(U \otimes k^\alg) \to \pi_1(U) \to \Gal(k^\alg/k) \to 1,\]
which we abbreviate by $\pi_1(U/k)$. Due to neglecting  base points, the extension $\pi_1(U/k)$ is only functorial in $U/k$ if regarded as an extension of $\Gal_k = \Gal(k^\alg/k)$ with $\pi_1(U \otimes k^\alg)$-conjugacy classes of maps between extensions.  

A $k$-rational point $u \in U(k)$ yields by functoriality a section $s_u$ of $\pi_1(U/k)$, with image the decomposition group of a point $\tilde{u}$ above $u$ in the universal pro-\'etale cover of $U$. Having neglected base points and due to the choice of $\tilde{u}$, only the class of a section up to conjugation by elements from $\pi_1(U \otimes k^\alg)$ is well defined. Let us denote by $S_{\pi_1{(U/k)}}$ the set of conjugacy classes of sections of $\pi_1(U/k)$. The section conjecture of Grothendieck's anabelian geometry \cite{letter}  speculates the following.
\begin{conj}[Grothendieck]
The natural map $U(k) \to S_{\pi_1(U/k)}$ which sends a rational point to the section given by its conjugacy class of decomposition groups is bijective if $U/k$ is a geometrically connected, smooth, projective curve of genus at least $2$ over an algebraic number field. % $k$.
\end{conj} 

This conjecture is wide open in general, although analogues  have been obtained over finite fields by Tamagawa \cite{Ta}, and over non constant $p$-adic function fields by Mochizuki \cite{Mz1}, as well as birational variants by Koenigsmann \cite{Koenigsmann} and more recently Pop \cite{P}. An approach using the Tannakian formalism has been advocated recently by Esnault and Hai \cite{eh}.

The injectivity part in the section conjecture was known to Grothendieck and consists only of an exercise exploiting the Mordell-Weil Theorem, see 
\cite{Mz1} Theorem 19.1, \cite{eh} Theorem 5.1, or Appendix \ref{app:injsec} Theorem \ref{thm:injsec}. 

\subsection{}
From now on we require \textbf{$k$ to be a field of characteristic $0$} and work with the full pro-finite fundamental groups only. Establishing tame or logarithmic as well as pro-$\ell$ versions of the material presented here involves only trivial modifications.

\subsection{Tangential and cuspidal sections}  Let $U$ be the complement of a divisor $Y$ with normal crossing in a geometrically connected, smooth proper variety $X$.  The fundamental group of the connected scheme of nearby points $U_y = U \times_X \Spec(\OO^{\rh}_{X,y})$ at a geometric point $y \in Y(k^\alg)$  can be computed by Zariski-Nagata purity and Abhyankar's Lemma to be the extension
\[ 1 \to \Hom\big(\OO^\ast(U_y)/(\OO_{X,y}^{\rh})^\ast, \hat{\bZ}(1)\big)  \to \pi_1(U_y) \to \Gal_{\kappa(y)} \to 1,\]
which we abbreviate by $\pi_1(U_y/\kappa(y))$, and where $\kappa(y) \subset k^\alg$ is the residue field of the closed point in $Y$ underlying $y$.
The inclusion $U_y \to U$ yields a natural map of extensions 
\[ \pi_1(U_y/\kappa(y)) \to \pi_1(U/k).\]

A \textbf{cuspidal} section of $\pi_1(U/k)$  is a section that factors over some $\pi_1(U_y)$ for a point $y \in Y$ with residue field $\kappa(y) =k$. Being cuspidal does not depend on the representative chosen in a given conjugacy class of sections. There is also a version of the section conjecture for affine curves stating that all but cuspidal sections come from rational points.
The condition on the genus gets replaced by asking the Euler-characteristic to be negative. 

In \cite{deligne} \S15.13-26 Deligne describes \textbf{tangential} base points for fundamental groups.  Taking the associated graded local rings and principal froms of maps describes an exact functor from finite \'etale covers of $U$ to finite covers of the tangent space $\rT_yX$ which may only branch along the tangent cone ${\rm TC}_yY$ of the normal crossing divisor $Y$. Hence we obtain a fibre above a tangent vector $v \in \rT_yX$ which is not tangent to a component of $Y$. For $k$-rational $v$ this construction leads to a section of $\pi_1(U/k)$, %see Appendix \ref{app:ts}, 
which clearly is cuspidal. Sections coming from suitable tangential base points will be called \textbf{tangential} sections.

\subsection{A guide through the paper}
This paper deals with the anabelian geometry of cuspidal sections in the case of curves with an emphasis on imposing only mild assumptions on the base field. 

There are three main topics. In Sections \ref{sec:kummer}--\ref{sec:abundance} we deal with
the question raised and firmly claimed by Grothendieck in \cite{letter}  of how many inequivalent cuspidal sections there are per boundary point. Contrary to Esnault and Hai in \cite{eh} where Tannakian methods are applied, we give a fairly elementary treatment\footnote{The author is indebted to Tam\'as Szamuely for inquiring such a 'pro-finite' presentation.} over fairly general base fields. 

In Sections \ref{sec:naka}--\ref{sec:units} we pursue a general investigation on additional structure of a fundamental group of a curve defining anabelian versions of a weight filtration\footnote{The idea of an anabelian weight filtration is due to Hiroaki Nakamura \cite{naka4} \S2.1.}, an orientation and units. In some sense these sections consist of the distillate of the author's careful reading of the mentioned papers by Nakamura \cite{naka1},  \cite{naka2},  \cite{naka3},  \cite{naka4}.

As an application of the general investigation sketched above, we obtain in Sections \ref{sec:genus0again} and \ref{sec:moduli} applications to weak anabelian geometry where the fundamental group shall determine the isomorphy type of a curve. The particular feature of our result especially in Section \ref{sec:moduli} is that we show directly how the moduli parameters in question are encoded in the arithmetic structure of the fundamental group extension.

What is unfortunately missing from the discussion so far, although the author would love to include these aspects into the theory, are the following topics:  (1) beyond curves when $\dim > 1$, (2) recognize tangential sections within cuspidal sections, and (3) finiteness and integrality of the anabelian degree and ramification indices  for arbitrary open maps.

We close the introduction by fixing some terminology as follows.

\subsection{} The \textbf{geometric fundamental group} of a geometrically connected variety $X/k$ is the group $\pig_1(X) =\pi_1(X \otimes k^\alg)$.

\subsection{} A \textbf{number field} $k$ is an algebraic extension of $\bQ$, whereas a finite extension of $\bQ$ will be called an \textbf{algebraic number field}.

\subsection{Neighbourhoods} A \textbf{neighbourhood} of a section $s \in S_{\pi_1(X/k)}$ is an open subgroup $H$ of $\pi_1(X)$ together with a representative of $s$ whose image is contained in $H$ considered up to conjugation by $H \cap \pig_1(X)$. Equivalently, a neighbourhood is a connected finite \'etale map $f:X' \to X$ and a section $s' \in S_{\pi_1(X'/k)}$ which maps to $s$ under $\pi_1(f)$.

%%%%%%%%%%%%%%%%%%%%%%%%%%%%%%%%%%%%%%%%%%%

\section{Fields with nontrivial Kummer theory} \label{sec:kummer}

\subsection{Pro-$\bN$ completion}
For an abelian group $A$ we set $\widehat{A} = \varprojlim_\bN A/nA$ and call it the pro-$\bN$ completion of $A$. If $A$ is finitely generated as a $\bZ$ module then $\widehat{A} = A \otimes \hat{\bZ}$.

\subsection{Definition}
We say that a field $k$ has \textbf{nontrivial Kummer theory}  if 
\[\widehat{k^\ast} = \rH^1\big(k, \hat{\bZ}(1)\big)\]
is infinite. Here $\rH^1\big(k, \hat{\bZ}(1)\big)$ is continuous cohomology of pro-finite groups with continuous cochains or, what amounts to the same, cohomology in the sense of Jannsen, see  \cite{Jannsen} Theorem 2.2.
 The isomorphism with $\widehat{k^\ast}$ is given by Kummer theory.

Among fields with trivial Kummer theory are finite fields, real closed fields and the maximal solvable extension $k^\solv$ of an arbitrary field $k$, for example $\bQ^\solv$.
\begin{lem}
Algebraic number fields, local fields and more generally fields $k$ that admit a discrete valuation $\nu$ with nondivisible value group $\Gamma$ have nontrivial Kummer theory.
\end{lem}
\begin{pro}
The valuation extends to a surjective map $\hat{\nu} : \widehat{k^\ast} \to \widehat{\Gamma}$ with infinite    image.
\end{pro}

\subsection{Characterization} The notation for the  torsion subgroup of a abelian group $A$ is $A_\tors$.
\begin{lem} \label{lem:finlevel}
Let $k$ be a field. Then $\rH^1(k,\bQ/\bZ(1)) = 0$ if and only if the natural map $\widehat{k^\ast_\tors} \to \widehat{k^\ast}$ is an isomorphism. These properties hold in particular when the natural projection $\widehat{k^\ast} \to k^\ast/(k^\ast)^n$ is an isomorphism for some $n \in \bN$.
\end{lem}
\begin{pro}
From the comparison of the multiplication by $n$ sequence of $\bQ/\bZ(1)$ and $\Gm$ we conclude that $k^\ast_\tors/(k^\ast_\tors)^n = k^\ast/(k^\ast)^n$ if and only if $\rH^1(k,\bQ/\bZ(1))$ has no $n$-torsion. This proves the claimed equivalence.

Let $\widehat{k^\ast} \to k^\ast/(k^\ast)^n$ be an isomorphism. Then for every $a \in k^\ast$ and $m\in \bN$ we have  $\zeta \in \boldsymbol{\mu}_n(k)$ such that $a\zeta \in (k^\ast)^m$. As the group $\mu_n(k)$ is finite, one $\zeta$ is good for all $m$, so $a \equiv \zeta^{-1}$ in $\widehat{k^\ast}$ and the anyway injective map $\widehat{k^\ast_\tors} \to \widehat{k^\ast}$ is also surjective.
\end{pro}
\begin{thm} \label{thm:2}
Let $k$ be a field. Then the following are equivalent.
\begin{enumerate}
\item[(a)] $\widehat{k^\ast}$ is finite, cyclic and generated by a root of unity,
\item[(b)] $\widehat{k^\ast}$ is finite,
\item[(c)] $\widehat{k^\ast}$ is countable,
\item[(d)] $\rH^1(k,\bQ/\bZ(1)) = 0$ and $k^\ast_\tors= \rH^0(k,\bQ/\bZ(1))$ is finite modulo its maximal divisible subgroup.
\end{enumerate}
In particular $k$ has nontrivial Kummer theory if and only if $\widehat{k^\ast}$ is uncountable.
\end{thm}
\begin{pro}
Condition (a) and (d) are equivalent by Lemma \ref{lem:finlevel}. 
Clearly (a) implies (b) implies (c). If (c) holds, the projective system defining
$\widehat{k^\ast}$ must stabilize at finite level, say $\widehat{k^\ast} \cong k^\ast/(k^\ast)^n$. From Lemma \ref{lem:finlevel} we deduce that $\widehat{k^\ast}$ is generated by roots of unity and isomorphic to $k^\ast_\tors/(k^\ast_\tors)^n$ which is a finite group, so (a) holds.
\end{pro}
\begin{cor}
The field $k$ has nontrivial Kummer theory if there is an $a \in k$ and $n\in \bN$ such that $k(\sqrt[n]{a})/k$ is not an abelian field extension.
\end{cor}
\begin{pro}
For such an $a \in k$ the Galois group of $\bigcup_m k(\sqrt[m]{a},\mu_m)$ over $k$ is nonabelian, hence the extension is not contained in the cyclotomic extension and $a$ is not contained in the kernel of $\widehat{k^\ast} \to \widehat{k(\mu_\infty)^\ast}$, which contadicts (a) of Theorem \ref{thm:2}.
\end{pro}

\subsection{Change of field} 
\begin{prop} \label{prop:cof}
Let $k$ be a field. The following are equivalent:
\begin{itemize}
\item[(a)] The extension $k(\mu_{2^\infty})/k$ is infinite or $\mu_4$ is contained in $k$.
\item[(b)] For any finite extension $E/F$ of fields which are finite algebraic over $k$ the natrual map $\widehat{F^\ast} \to \widehat{E^\ast}$ is injective.
\end{itemize}
\end{prop}
\begin{pro}
In (b) we may assume that $E/F$ is Galois. By the inflation/restriction sequence the kernel of $\widehat{F^\ast} \to \widehat{E^\ast}$ equals $\rH^1(E/F, \hat{\bZ}(1)\big(E))$. We may restrict to the $p$-part, and after the usual corestriction argument with the $p$-Sylow subgroup and d\'evissage we see that in (b) it is enough to ask $\rH^1(E/F,{\bZ}_p(1)\big(E))=0$ for all $p$-cyclic Galois extensions $E/F$.

The group of coefficients ${\bZ}_p(1)\big(E)$ either vanishes or $\mu_{p^\infty}$ is contained in $E$.
For $p$ odd, either way the action of $\Gal(E/F)$ must be trivial as ${\bZ}_p^\ast$ contains no $p$-torsion. Hence in this case
\[\rH^1(E/F, {\bZ}_p(1)\big(E)) = \Hom(\Gal(E/F), {\bZ}_p(1)\big(E)) = 0.\]
For $p=2$ we may argue analogously except in the case that $\Gal(E/F)$ acts via complex conjugation, which means $\mu_4 \not\subset F$ and $\mu_{2^\infty} \subset E$. Therefore (a) implies (b).

But if (a) fails, then there is a quadratic extension $E/F$ finite over $k$ with $\mu_{2^\infty} \subset E$ and $E=F(i)$ such that
\[ \ker\big(\widehat{F^\ast} \to \widehat{E^\ast}\big) = \rH^1(\bZ/2\bZ, \hat{\bZ}(1)) = \ker\big(\widehat{\bR^\ast} \to \widehat{\bC^\ast}\big) = \{\pm 1\}\]
contradicting (b).
\end{pro}
\begin{cor} \label{cor:cof}
For a finite field extension $E/F$ the kernel of $\widehat{F^\ast} \to \widehat{E^\ast}$ has order at most $2$ with the class of $-1$ being the only possibly nontrivial element.
\end{cor}
\begin{pro}
By Proposition \ref{prop:cof} the kernel is contained in the corresponding kernel for $E=F(i)$, which at most equals the kernel for the extension $\bC/\bR$ as shown in the proof above. 
\end{pro}

%%%%%%%%%%%%%%%%%%%%%%%%%%%%%%%%%%%%%%%%%%%

\section{Inertia subgroups} \label{sec:inertia}
\subsection{Cusps}
From now on $U$ will be the complement of a reduced divisor $Y$ in a geometrically connected, smooth, projective curve $X/k$.  The Euler-characteristic of $U $ is  $\chi(U) = 2 - 2g - \deg(Y)$, where $g$ is the genus of $X$. 

We fix a pro-universal cover $\tilde{U} \to U \otimes_k k^\alg $ together with an isomorphism of $\Gal_k$-extensions of $\pi_1(U/k)$ with the opposite groups of 
\[ 1 \to \Aut(\tilde{U}/U \otimes_k k^\alg) \to \Aut(\tilde{U}/U) \xrightarrow{\pr} \Aut(\Spec(k^\alg)/\Spec(k)) \to 1.\]
The normalisation of $Y$ in $\tilde{U}$ is the $k^\alg$-scheme $\tilde{Y}$.
The set of \textbf{cusps} of $U$ is $\cp(U) = \Hom_k(\Spec(k^\alg),Y)$  as a set with  $\Gal_k$ action. The set of prolongations of cusps to $\tilde{U}$ is $\tcp(U) = \Hom_{k^\alg}(\Spec(k^\alg), \tilde{Y})$,  which carries a natural $\pi_1(U)$ action by 
\[ \gamma . \tilde{y} :=  \gamma^{-1} \circ \tilde{y} \circ \pr(\gamma)\]
The projection map $\tcp(U) \surj \cp(U), \tilde{y} \mapsto y$ is equivariant and the quotient map for the induced  $\pig_1(U)$ action. 
The inertia group $\rI_\ty$ (resp.\ the decomposition group $\rD_\ty$) of a cusp $\tilde{y} \in \tcp(U) $ is the stabilizer under the action of $\pig_1(U)$ (resp.\ $\pi_1(U)$). Hence for $\gamma \in \pi_1(U)$ we have 
\[ \rI_{\gamma\tilde{y}/\gamma y} = \gamma \rI_\ty \gamma^{-1} \qquad \text{and} \qquad \rD_{\gamma\tilde{y}/\gamma y} = \gamma \rD_\ty \gamma^{-1}.
\]
Moreover, the choice of $\tilde{y} \in \tcp(U)$ determines a representative for the map of extensions $\pi_1(U_y/\kappa(y)) \to \pi_1(U/k)$ with image the extension $\cD_\ty$ given by
\[ 1 \to \rI_{\tilde{y}/y} \to \rD_{\tilde{y}/y} \to \Gal_{\kappa(y)} \to 1.\]
\begin{prop} \label{prop:localsplit}
The extension $\pi_1(U_y/\kappa(y))$ splits and there is a natural free transitive action on $S_{\pi_1(U_y/\kappa(y))}$ by $\widehat{\kappa(y)^\ast}$ .
\end{prop}
\begin{pro}
Fix a parameter $t$ at $y$. Adjoining the $n^{th}$ root of $t$ defines a compatible family of finite \'etale covers of $U_y$ and correspondingly a family of open subgroups of $\pi_1(U_y)$ the intersection of which defines a splitting.
By Proposition \ref{prop:h1} of Appendix  \ref{app:h1} the set $S_{\pi_1(U_y/\kappa(y))}$ carries a free transitive action by $\rH^1(\kappa(y),\hat{\bZ}(1)) = \widehat{\kappa(y)^\ast}$.
\end{pro}

%---------------------------------------------------------------------------------------------------------------------------------

\subsection{Universal ramification}
From the topological computation of fundamental groups via GAGA we get the following exact sequence of $\Gal_k$ modules
\begin{equation} \label{eq:piab}
%0 \to 
\hat{\bZ}(1) \xrightarrow{\Delta} \hat{\bZ}(1) \otimes \bZ[\cp(U)] \to \pig_1^\ab(U) \to \pig_1^\ab(X) \to 0,
\end{equation}
where moreover $\hat{\bZ}(1)\cdot y$ for a cusp $y \in \cp(U)$ maps to the image of $\rI_\ty$ in 
$\pig_1^\ab(U)$ regardless of the prolongation $\tilde{y}$ of $y$, and $\Delta$ is the diagonal map. The group $\pig_1^\ab(X)$ is a free $\hat{\bZ}$ module of rank $2g$.
\begin{lem} \label{lem:cusptoinfty}
If $\chi(U) < 0$ and $U$ is affine, then $\# \cp(V)$ tends to infinity as $V$ ranges over connected finite \'etale covers of $U$ that, moreover,  we may choose among the neighbourhoods of a  section $s\in S_{\pi_1(U/k)}$.
\end{lem}
\begin{pro}
Otherwise $\rI_\ty$ had finite index in $\pig_1(U)$ which contradicts (\ref{eq:piab}).
\end{pro}
\begin{lem} If  $\chi(U) \leq 0$ then $\pi_1(U_y/\kappa(y)) \cong \cD_\ty$ injects into $\pi_1(U/k)$ for each cusp $y$ of $U$.
\end{lem}
\begin{pro}
As $\hat{\bZ}(1)$ is torsion free we may replace $U$ by a finite \'etale cover. The lemma follows from (\ref{eq:piab}) as soon as $\#\cp(U) \geq 2$, which leaves the trivial case $U \otimes k^\alg \cong \Gm$.
\end{pro}

%---------------------------------------------------------------------------------------------------------------------------------

\subsection{Local theory}
From now on we assume that $U$ is a hyperbolic curve, i.e., that $\chi(U)$ is negative.
\begin{prop} Let $\rI_{\tilde{y}_1/y_1} \cap \rI_{\tilde{y}_2/y_2} \not= 1$ for some cusps $\tilde{y}_1, \tilde{y_2} \in \tcp(U)$.  Then already $\tilde{y_1} = \tilde{y}_2$.
\end{prop}
\begin{pro}
Arguing by contradiction we may replace $U$ by a finite \'etale cover and assume $y_1 \not=  y_2$ with $\# \cp(U) \geq 3$. Then  $\rI_{\tilde{y}_1/y_1}$ and $\rI_{\tilde{y}_2/y_2}$ intersect trivially in $\pig_1^\ab(U)$  by (\ref{eq:piab}).
\end{pro}
\begin{cor}  \label{cor:normal}
The normalizer of $\rI_\ty$ in $\pi_1(U \otimes k^\alg)$  (resp.\ $\pi_1(U)$) are 
$\rI_\ty$ (resp.\ $\rD_\ty$).
\end{cor}
\begin{pro}
If $\gamma \in \pi_1(U)$ normalises $\rI_\ty$, then $\rI_\ty \cap \rI_{\gamma.\tilde{y}/\gamma.y} \not= 1$, hence $\gamma.\tilde{y} = \tilde{y}$ and thus $\gamma \in \rD_\ty$. Furthermore, we have $\rI_\ty = \rD_\ty \cap \pig_1(U)$.
\end{pro}
\begin{cor} \label{cor:abel}
Let $A$ be an abelian subgroup of $\pig_1(U)$ that nontrivially intersects an inertia subgroup $\rI_\ty$. Then $A$ is contained in $\rI_\ty$.
\end{cor}
\begin{pro}
Let $\gamma$ be a nontrivial element of $A \cap \rI_\ty$ and $a \in A$ arbitrary. Then the intersection of $\rI_\ty$ with $a\rI_\ty a^{-1} = \rI_{a.\tilde{y}/a.y}$ contains $\gamma$, hence $a.\tilde{y} = \tilde{y}$ and $a \in \rI_\ty$.
\end{pro}

%%%%%%%%%%%%%%%%%%%%%%%%%%%%%%%%%%%%%%%%%%%

\section{Abundance of cuspidal sections}  \label{sec:abundance}

\subsection{The anabelian proof} In this section we give two proofs for the fact that quite generally the set of cuspidal sections is uncountable, a fact emphasized by Esnault and Hai in \cite{eh}  Corollary 6.9.
\begin{thm} \label{thm:1}
Let us assume that the injectivity part of the section conjecture holds for proper curves of genus at least $2$ over finite extensions of $k$. Then the following holds.
\begin{enumerate}
\item A cuspidal section uniquely determines the cusp $y$ for which it factors over $\pi_1(U_y)$.
\item Let $y$ be a $k$-rational cusp of $U$. The natural map $S_{\pi_1(U_y/k)} \to S_{\pi_1(U/k)}$ is injective and the image disjoint from the sections induced by rational points $u \in U(k)$ or cuspidal sections originating from a different cusp.
\end{enumerate}
\end{thm}
\begin{rmk} (1)
The image of $S_{\pi_1(U_y/k)} \to S_{\pi_1(U/k)}$ is called the \textbf{`Paket' (engl.\ packet)} associated to the cusp $y$ in \cite{letter}.

(2) Let $k$ be a field such that for any abelian variety $A/k$ the group of $k$-rational poins $A(k)$ has no divisible elements. Then the assumption in Theorem \ref{thm:1}  holds, see Appendix \ref{app:injsec}.
\end{rmk}
\begin{pro}
Let $s$ be a section. Replacing $U$ by a neighbourhood of $s$ we may assume that $X$ has genus at least $2$. 

(1) If the image of $s$ is contained in $\rD_{\tilde{y}_1/y_1} $ and  $\rD_{\tilde{y}_2/y_2}$ then $s$ maps simultaneously to $s_{y_1}$ and $s_{y_2}$ under $S_{\pi_1(U/k)} \to S_{\pi_1(X/k)}$, hence $y_1=y_2$ by assumption.

(2) By the argument above it remains to show that sections $s,t$ of $\pi_1(U_y/k)$ which become conjugate as sections of $\pi_1(U/k)$ are already conjugate in $\pi_1(U_y/k)$. Let $\gamma \in \pig_1(U)$ be such that $s=\gamma()\gamma^{-1} \circ t$. Then for suitable $\tilde{y}$ the image $s(Gal_k)$ is contained in $\rD_\ty$ and $\gamma(\rD_\ty)\gamma^{-1} = \rD_{\gamma.\tilde{y}/\gamma.y}$. As above, we conclude that $y=\gamma.y$. Applying the same reasoning to all finite \'etale covers of $U$ we also deduce $\gamma.\tilde{y} = \tilde{y}$. Hence $\gamma \in \rI_\ty$ and so $s$ and $t$ are already conjugate as sections of $\pi_1(U_y/k)$.
\end{pro}

\begin{cor} If the injectivity part of the section conjecture holds for proper curves of genus at least $2$ over finite extensions of $k$ and if $k$ has a nontrivial Kummer theory then 
 each $k$-rational cusp gives rise to an uncountable set of cuspidal sections.
\end{cor}
\begin{pro}
This follows from Theorem \ref{thm:1} in conjunction with Proposition \ref{prop:localsplit}.
\end{pro}

%---------------------------------------------------------------------------------------------------------------------------------

\subsection{The abelian proof}
The abelianised fundamental group extension $\pi_1^\ab(U/k)$ is the push\-out of the extension $\pi_1(U/k)$ by the characteristic quotient $\pig_1(U) \surj \pig_1^\ab(U)$. 
\begin{thm} \label{thm:3}
Let $k$ be a field with nontrivial Kummer theory such that one of the following two assumptions hold.
\begin{itemize}
\item[(i)] For any abelian variety $A/k$ the group of $k$-rational torsion points $A(k)_\tors$ contains no nontrivial divisible element.
\item[(ii)] The $\hat{\bZ}$ rank of $\widehat{k^\ast}$ is infinite.
\end{itemize}
Let $y$ be a $k$-rational cusp of $U$. The natural map $S_{\pi_1(U_y/k)} \to S_{\pi_1(U/k)}$ has uncountable image.
\end{thm}
\begin{rmk} \label{rmk:Qab}
Assumption (ii) holds for $k=\bQ^\ab$. For each prime $p$ we have a valuation $\nu_p$ on $\bQ^\ab$ with values in $1/(p-1)\bZ[1/p]$. For different $p$ these are mutually independent, whence for finite sets $S$ of primes a surjective map $\widehat{\bQ^{\ab,\ast}} \to \prod_{p \in S} \widehat{1/(p-1)\bZ[1/p]}$. Choosing a prime $\ell \not\in S$ we conclude that the $\bZ_\ell$ rank of $\widehat{\bQ^{\ab,\ast}} \otimes \bZ_\ell$ is at least $\# S$, and (ii) follows. 
\end{rmk}
\begin{pro}
Replacing $U$ by a neighbourhood of a cuspidal section at $y$ we may assume by Lemma 
\ref{lem:cusptoinfty} that $U$ has at least $2$ cusps. It is moreover sufficient to prove that the image of 
$S_{\pi_1(U_y/k)} \to S_{\pi^\ab_1(U/k)}$ is uncountable, which by Proposition \ref{prop:h1} of Appendix \ref{app:h1} amounts to the map $\widehat{k^\ast}=\rH^1(k,I_\ty) \to \rH^1(k,\pig_1^\ab(U))$ having uncountable image.

The $k$-rational cusp $y$ can be used to split the diagonal map $\Delta$ in \ref{eq:piab}, so that Galois cohomology yields the exact sequence
\[ \rH^0(k,\pig_1^\ab(X)) \xrightarrow{\delta} \big(\bigoplus_{z} \widehat{\kappa(z)^\ast}\big)/\widehat{k^\ast} \to \rH^1(k,\pig_1^\ab(U)),\]
where $z$ runs over a set of representatives of $\Gal_k$ orbits on $\cp(U)$. By Corollary \ref{cor:cof} the kernel of the map $\widehat{k^\ast} \to \big(\bigoplus_{z} \widehat{\kappa(z)^\ast}\big)/\widehat{k^\ast}$ onto the component of the cusp $y$ is at most  of order $2$. The result follows because in case (i) the map $\delta$ is trivial and in case (ii) it has image of finite rank, whereas $\widehat{k^\ast}$ has infinite rank.
\end{pro}

%%%%%%%%%%%%%%%%%%%%%%%%%%%%%%%%%%%%%%%%%%%

\section{Anabelian weight filtration following Nakamura} \label{sec:naka}

We recall the theory of the anabelian weight filtration after Nakamura \cite{naka2} \S3, \cite{naka4} \S2.1.
As before, $U$ will be the complement of a reduced divisor $Y$ in a geometrically connected, smooth, projective curve $X/k$.

\subsection{Definition} The \textbf{anabelian weight filtration} of $U$ is the subset $W_{-2}(U)$ of all inertia elements
\[ W_{-2}(U) = \bigcup_{\tilde{y} \in \tcp(U)} \rI_\ty \subset \pig_1(U)\]
The set $W_{-2}(U)$ is pro-finite, hence compact, and preserved under conjugation by $\pi_1(U)$.
Strictly speaking, the designation  \textit{anabelian} is only justified in instances where $W_{-2}(U)$ can be described in terms of the fundamental group alone.

\subsection{Definition}
We say that a pro-cyclic closed subgroup $I \subset \pig_1(U)$ is 
\textbf{essentially cyclotomically normalized} if there is a subgroup $N$ of the normalizer $N_{\pi_1(U)}(I)$ of $I$ in $\pi_1(U)$ which projects to an open subgroup in $\Gal_k$ and such that the induced action of $N$ on $I$ is via the cyclotomic character of $\Gal_k$. The pro-cyclic closed subgroup $I$ is called 
\textbf{cyclotomically normalized} if it is essentially cyclotomically normalized with $N=N_{\pi_1(U)}(I)$.

\subsection{Characterisation}
We say that a field $k$ \textbf{distinguishes tori from abelian varieties} if for all tori $S/k$ and all abelian varieties $A/k$ the group $\Hom_{\Gal_k}(\rT S,\rT A)$ of maps between Tate modules is trivial. Using Weil restriction of scalars, we see that if $k$ distinguishes tori from abelian varieties then so does any finite extension.  An equivalent definition occuring in \cite{CT} \S2 Definition, is the notion of the cyclotomic character being a non-Tate character.

Due to the theory of weights \`a la Deligne, algebraic number fields distinguish tori from abelian varieties. More precisely\footnote{The author thanks Akio Tamagawa for the reference to Ribet's Theorem.}, 
the following lemma is a direct reformulation of Theorem 1 of Ribet's appendix to \cite{KL}.

\begin{lem} \label{lem:Ribet}
A number field $k$ such that $k\bQ^\ab/\bQ^\ab$ is finite distinguishes tori from abelian varieties.
\end{lem}

\begin{thm}[Nakamura] \label{thm:naka}
Let $k$ be a field that distinguishes tori from abelian varieties, and let $U/k$ be a hyperbolic curve. Then the following are equivalent for a  nontrivial closed subgroup $\rI$ in $\pig_1(U)$:
\begin{enumerate}
\item[(a)] $\rI$ is contained in an inertia subgroup $\rI_\ty$,
\item[(b)] $\rI$ is pro-cyclic and cyclotomically normalized,
\item[(c)] $\rI$ is pro-cyclic and essentially cyclotomically normalized.
\end{enumerate}
\end{thm}
\begin{pro}
(a) implies (b): As $\rI_\ty$ is pro-cyclic, $\rI$ is pro-cyclic as well. If $\gamma \in \pi_1(U)$ normalizes $\rI$ then $\rI_{\gamma.\tilde{y}/\gamma.y} = \gamma \rI_\ty \gamma^{-1}$ intersects $\rI_\ty$ nontrivially in $\rI$. Hence $\gamma.\tilde{y} = \tilde{y}$ and $\gamma \in \rD_\ty =N_{\pi_1(U)}(\rI_\ty)$. As $\rI_\ty$ is cyclotomically normalized the same follows for $\rI$.

(b) implies (c) is trivial, hence it remains to conclude (a) from (c). We argue by contradiction. By Corollary \ref{cor:abel} we may replace $\rI$ by $\rI \otimes \bZ_\ell \cong \bZ_\ell$ for some prime $\ell$. Again by Corollary \ref{cor:abel} we have $I \cap W_{-2}(U) = 1$. Choose a compact set $C \subset \rI$ avoiding $1$. The Hausdorff property yields a normal subgroup $N$ of $\pig_1(U)$ such that $CN/N$ has empty intersection with $W_{-2}(U)N/N$ in $\pig_1(U)/N$.  

The group generated by $W_{-2}(U) \cap N\rI$ in $N\rI/N$ is cyclic of order a power of $\ell$ hence the image of some $\rI_\ty \cap N\rI$ which therefore avoids $CN/N$. Thus the intermediate $\ell$-cyclic cover $U''\to U'$ corresponding to $N(\rI)^\ell \subset N\rI$ is unramified over the corresponding cusps $X'-U'$. For some finite field extension we obtain a nontrivial Galois invariant map from $\rI\cong \bZ_\ell(1) \to \rT_\ell A'$ where $A'$ is the Albanese variety of $X'$. This contradicts the assumption that $k$ distinguishes tori from abelian varieties.
\end{pro}

We deduce an anabelian description of inertia subgroups, hence the set $\tcp(U)$ and the anabelian weight filtration $W_{-2}(U)$.
\begin{cor} \label{cor:anWF}
Let $k$ be a field that distinguishes tori from abelian varieties, and let $U/k$ be a hyperbolic curve. Then the set of inertia subgroups in $\pig_1(U)$ coincides with the set of maximal closed subgroups of $\pig_1(U)$ among those which are pro-cyclic and cyclotomically normalized.
\end{cor}
\begin{pro}
An inertia subgroup $\rI_\ty$ is pro-cyclic and cyclotomically normalized, and in fact a maximal such subgroup as otherwise $\rI_\ty $ were properly contained in an abelian subgroup contradicting Corollary \ref{cor:abel}. 

Conversely, a maximal pro-cyclic and cyclotomically normalized subgroup $\rI$ is contained in some $\rI_\ty$ by Theorem \ref{thm:naka} and must in fact coincide with $\rI_\ty$ because of maximality.
\end{pro}
\begin{cor}
Let $k$ be a field that distinguishes tori from abelian varieties, and let $U/k$ be a hyperbolic curve. A section $s \in S_{\pi_1(U/k)}$ is cuspidal if and only if the image $s(\Gal_k)$ cyclotomically normalizes a pro-cyclic subgroup of $\pig_1(U)$.
\end{cor}
\begin{pro}
A cuspidal section at the cusp $y$ cyclotomically normalizes $\rI_\ty$ for some choice of $\tilde{y}$ depending on various base points and paths. 

Conversely, let the image of the section $s$ cyclotomically normalizes the pro-cyclic subgroup $\rI$, which then by Theorem \ref{thm:naka} must be contained in some inertia subgroup $\rI_\ty$.
As above, we see that $s(\Gal_k)$ also normalizes $\rI_\ty$. Thus $s$ factors over $\rD_\ty=\pi_1(U_y)$ and is cuspidal.
\end{pro}

\subsection{Weight preserving maps} 
A \textbf{weight preserving} map is a continuous group homomorphism
$\ph:\pi_1(U/k) \to \pi_1(V/k)$ of fundamental groups of hyperbolic curves $U$ (resp.\ $V$) over $k$ which are complements of divisors $Y$  (resp.\ $Z$) on geometrically connected, smooth projective curves $X$  (resp.\ $W$), such that $\ph$ restricted to $W_{-2}(U)$ takes values in $W_{-2}(V)$.
This does not exclude that an inertia group of $\pi_1(U)$ may lie in the kernel of $\ph$. 
From  Corollary \ref{cor:abel} we learn that $\ph$ induces a partially defined $\pi_1$-equivariant map $\tcp(V) \to \tcp(V)$, which is defined precisely on those cusps $\ty$ for which $\rI_\ty$ is not killed by $\ph$.

The goal of local theory in anabelian geometry is among others to show that any group homomorphism automatically preserves weights.  For hyperbolic curves over fields that distinguish between tori and abelian varieties this was achieved by Nakamura as recalled above. 
As we have different applications in mind we will assume in subsequent sections that our group homomorphisms preserve weights in the first place.

%%%%%%%%%%%%%%%%%%%%%%%%%%%%%%%%%%%%%%%%%%%

\section{Anabelian local cohomology for curves} \label{sec:loccoh}

%\subsection{Modified Hochschild-Serre spectral sequence} \label{subsec:modsf}
%Let $G$ be a profinite group with closed subgroup $I$, normal subgroup $N$ generated by $I$ and quotient $\Gamma = G/N$. The functor of $N$-invariants on discrete $G$-modules is identical to the functor of $I$-invariants. In the Hochschild-Serre spectral sequence 
%\[  \rH^p\big(\Gamma,\rH^q(N,M)\big) \Longrightarrow \rH^{p+q}(G,M)\]
%we may therefore replace $N$-cohomology by $I$-cohomology. If, moreover, $I$ has cohomological dimension $1$, the spectral sequence degenerates into a long exact sequence of cohomology
%\[  0  \to  \rH^1\big(\Gamma,M^I\big) \to  \rH^1\big(G,M\big) \to  \rH^0\big(\Gamma,\rH^1(I,M)\big) \xrightarrow{d_2^{0,1}} \rH^2\big(\Gamma,M^I\big) \to \rH^2\big(G,M\big)  \dots \]
%where $M^I$ are the invariants under $I$ of the coefficient $G$-module $M$.

%%%% the map $d_2^{0,1}$ equals the transgression, see [NSW2]. But I am unable to identify
%%%% this with the boundary operator from local cohomology sequence, at least not without
%%%% the general nonsense on comparison with mapping cone categories of sheaves.

\subsection{Local cohomology sequence}
We will give another  anabelian definition of the sequence (\ref{eq:piab}). The local cohomology sequence for the pair $(X \otimes k^\alg,U \otimes k^\alg)$ in \'etale cohomology with coefficients in $\bQ/\bZ$ has the following Pontrjagin dual
\begin{equation} \label{eq:loc}
 \rH^2(X \otimes k^\alg,\bQ/\bZ)^\vee \xrightarrow{\Delta} \bigoplus_{y \in \cp(U)} 
 \rH_y^2(X \otimes k^\alg,\bQ/\bZ)^\vee  \to \pig_1^\ab(U) \to \pig_1^\ab(X ) \to 0,
\end{equation}
which is even injective on the left except if $X=U$. Let us assume that $U$ is hyperbolic.

\subsection{Excision} Let $U_y^{\sh}$ be the scheme of geometric nearby points of $y$, which is the preimage of $U$ in $X_y^{\sh}  = \Spec(\OO_{X,y}^\sh)$. Naturality of (\ref{eq:loc}) for the map of pairs $(X_y^\sh,U_y^\sh) \to (X,U)$ and excision compute the local cohomology group 
$\rH_y^2(X \otimes k^\alg,\bQ/\bZ)^\vee = \rH_y^2(X_y^\sh,\bQ/\bZ)^\vee $
together with the map to $\pi_1^\ab(U \otimes k^\alg)$ as the natural map of an inertia group 
\[ \rI_\ty = \pi_1^\ab(U_y^\sh) \to \pig_1^\ab(U) \]
for a choice of cusp $\tilde{y}$ above $y$ reflected in the choice of base points.  The image is independent of the choice  of the prolongation $\tilde{y}$ and will be denoted by $\rI_y$.

\subsection{The orientation module}
 Apart from the map $\Delta$ all constituents of (\ref{eq:loc}) have now an anabelian definition in terms of the fundamental group $\pi_1(U/k)$ together with its weight filtration alone. Here is the anabelian definition of the map $\Delta$. 
 
For $X$ of genus $0$ we define the \textbf{orientation module} $\orient_{\pi_1(U/k)}$ and $\Delta$ through the exactness of the sequence 
 \[ 0 \to \orient_{\pi_1(U/k)}  \xrightarrow{\Delta} \bigoplus_{y \in \cp(U)} \rI_y \to \pig_1^\ab(U).\]
If  the genus of $X$ is at least $1$, then $X$ is an \'etale $K(\pi,1)$, cf. \cite{Stix} Appendix A, and we define the \textbf{orientation module}   as  
\[ \orient_{\pi_1(U/k)} = \rH^2(\pig_1(X),\bQ/\bZ)^\vee =  \rH^2(X \otimes k^\alg,\bQ/\bZ)^\vee .\]
The $y$ component of the map $\Delta$ will be defined by functoriality from the case of the pair $(X,X-\{y\})$ and as the Pontrjagin dual to a map
\[ \delta_y : \Hom(\rI_y, \bQ/\bZ)  \to \rH^2(\pig_1(X ),\bQ/\bZ) \]
as follows. There is a maximal intermediate quotient $\pig_1((X-\{y\}) ) \surj E \xrightarrow{\alpha} \pig_1(X)$ such that $\ker(\alpha)$ is central in $E$. Then canonically $\ker(\alpha)=I_y$ and we have a central extension 
\begin{equation} \label{eq:univcentral}
1 \to I_y \to E \to \pig_1(X) \to 1.
\end{equation}
Pushing $E$  by $\chi \in \Hom(\rI_y,\bQ/\bZ)$ defines the central extension $\delta_y(\chi) \in \rH^2(\pig_1(X),\bQ/\bZ) $.

\subsection{Comparison} \label{subsec:comp}
It remains to compare the anabelian definition of $\Delta$ with the geometric definition, which we will do
restricted to the $n$-torsion part and twisted by $\mu_n$. In the following diagram
\[
\xymatrix{
\rH^0\big(U_y^\sh,\Gm\big) \ar@{->>}[r]^{\delta_\loc} \ar@{->>}[d]^{\delta_\kum} \ar@{}[dr]|{\fbox{-1}} &
\rH_y^1\big(X_y^\sh,\Gm\big)  \ar[d]^{\delta_\kum} &  
\ar[l]_(0.55){\exc}^(0.55)\cong  \rH_y^1\big(X \otimes k^\alg,\Gm\big)  \ar[r]   \ar[d]^{\delta_\kum} &
\rH^1\big(X \otimes k^\alg,\Gm\big)  \ar[d]^{\delta_\kum} \\
\rH^1\big(U_y^\sh,\mu_n\big)  \ar[r]^{\delta_\loc}_{\cong} \ar@{}[drrr]|{\fbox{\qquad \qquad \qquad -1\qquad \qquad \qquad }} & 
\rH_y^2\big(X_y^\sh,\mu_n\big) & 
\ar[l]_(0.55){\exc}^(0.55)\cong  \rH_y^2\big(X \otimes k^\alg,\mu_n\big)  \ar[r]^{\Delta_y^\vee}   &  
\rH^2\big(X \otimes k^\alg,\mu_n\big) \\
\rH^1\big(\pi_1(U_y^\sh),\mu_n\big)  \ar@{=}[r] \ar[u]_\cong^{\comp}
& \Hom(I_y,\mu_n) \ar[rr]^{\delta_y = {\rm push \ } E} & & 
\rH^2\big(\pig_1(X), \mu_n\big)  \ar[u]_\cong^{\comp}
}
\]
the upper left facet commutes up to sign by \cite{sga4.5} cycle 2.1.3. We compute the maps going around the diagram along the outside border and establish that the bottom facet also commutes only up to a sign. Let $f$ be a parameter of $X$ at $y$. The mod $n$ tame character $\in  \Hom(I_y,\mu_n)$
maps under the comparison map to the $\mu_n$ torsor $[\sqrt[n]{f}] \in \rH^1\big(U_y^\sh,\mu_n\big) $ of $n^{th}$ roots of $f$, which equals $\delta_\kum(f)$ by \cite{sga4.5} cycle 1.1.1. By \cite{sga4.5} cycle 1.1.4, we may interpret $\delta_\loc(f) \in \rH_y^1\big(X_y^\sh,\Gm\big)$ as the trivial line bundle on $X_y^\sh$ plus the trivialisation given by $f$ over $U_y^\sh$ which maps under excision and the natural map to the class of  line bundle $\OO(y) \in \rH^1\big(X \otimes k^\alg,\Gm\big)$ associated to the divisor $y$. By the very definition of \cite{sga4.5} cycle 2.1.2, the class $\delta_\kum(\OO(y))$ is the cycle class $\cl_y$ of $y$ and the first chern class $\rc_1(\OO(y))$. By \cite{Mz2} Lemma 4.5, the corresponding central extension in $\rH^2(\pig_1(X),\mu_n)$ is given by the mod $n$ push of the homotopy sequence 
\begin{equation} \label{eq:fundcentral}
 1 \to \hat{\bZ}(1) \to \pi_1(\bL^\circ) \to \pig_1(X) \to 1
\end{equation}
of the complement $\bL^\circ$ of the zero section in $\bL=\underline{\Spec}_{X \otimes k^\alg} \big(\Sym^\bullet \OO(-y)\big)$ which is the geometric line bundle whose sections naturally coincide 
with sections of $\OO(y)$. 

The section $1$ above $X - \{y\}$ defines a map $(X-\{y\}) \otimes k^\alg \to \bL^\circ \to X \otimes k^\alg$ which leads to a map of extensions 
\[ 
\xymatrix{1\ar[r] & I_y \ar[r]  \ar[d]^\chi & E \ar[r]  \ar[d] & \pig_1(X)  \ar@{=}[d] \ar[r] & 1 \\
 1 \ar[r] & \hat{\bZ}(1) \ar[r] & \pi_1(\bL^\circ) \ar[r] & \pig_1(X) \ar[r] & 1}
\]
and all that remains is to identify the character $\chi$ as the tame character. The trivialisation of $\bL$ over $X_y^\sh$ via $f$ leads to a diagram 
\[
\xymatrix{
\Gm \ar@{=}[dr] \ar[r] & \bL^\circ|_{X_y^\sh} \ar[r]^(0.4)\cong \ar[d]^{\pr_1} & \Gm \times X_y^\sh  \\
  & \Gm  & U_y^\sh \ar[ul]_{1} \ar[l]_f  \ar[u]_{f,{\rm incl}} }
\]
which on fundamental groups gives
\[
\xymatrix{
\hat{\bZ}(1) \ar@{=}[dr] \ar[r]^(0.4)\cong & \pi_1(\bL^\circ|_{X_y^\sh}) \ar[r]^(0.4)\cong \ar[d]^{\pi_1(\pr_1)} & \pi_1(\Gm \times X_y^\sh)  \\
  & \hat{\bZ}(1)  & \pi_1(U_y^\sh) = I_y \ar[ul]_{\chi} \ar[l]_{\pi_1(f)}  \ar[u]_{\pi_1(f,{\rm incl})} }
\]
and reveils that $\chi = \pi_1(f)$ which is exactly the tame character.

% \[
%\xymatrix{\Hom(\rI_y,\mu_n) \ar[r]^(0.4){\delta_y} \ar[d]^\cong \ar@{}[dr]|{\fbox{-1}}  &  \rH^2(\pi_1(X \otimes k^\alg),\mu_n) \ar[d]^\cong \\
%\rH_y^2(X \otimes k^\alg,\mu_n)   \ar[r]^{\Delta^\vee}  &  \rH^2(X \otimes k^\alg,\mu_n) 
%}
%\]

%
% failed attempt: comparison with modified Hochschild-Serre spectral sequence
%
%as the Pontrjagin dual to the map
% \[ \Hom(\rI_y, \bQ/\bZ) = \rH^1(\pi_1(U_y^\sh),\bQ/\bZ) \to \rH^2(\pi_1(X \otimes k^\alg),\bQ/\bZ)\]
% that comes from the sequence in section \ref{subsec:modsf} above applied to $G = \pi_1(X \otimes k^\alg - \{y\})$, $\rI = \rI_\ty$ for some choice of prolongation $\tilde{y}$ and $\Gamma = \pi_1(X \otimes k^\alg)$ with $M = \bQ/\bZ$, the action of $\Gamma$ on $\rH^1(\pi_1(U_y^\sh),\bQ/\bZ)$ being trivial.
 
 \subsection{Anabelian local cohomology sequence}
 We conclude that the anabelian local cohomology sequence 
 \begin{equation} \label{eq:anabloc}
 \orient_{\pi_1(U/k)}  \xrightarrow{\Delta} \bigoplus_{y \in \cp(U)} 
 \rI_y \to \pig_1^\ab(U) \to \pig_1^\ab(X) \to 0,
\end{equation}
as constructed above is isomorphic to the local cohomology sequence in \'etale cohomology  (\ref{eq:loc}) up to a sign for the map $\Delta$. The map $\Delta$ is even injective except if $X=U$, and the groups $\orient_{\pi_1(U/k)}$ and $\rI_y$ are all isomorphic to $\hat{\bZ}(1)$ as Galois modules.

%%%%%%%%%%%%%%%%%%%%%%%%%%%%%%%%%%%%%%%%%%%

\section{Orientation and degree} \label{sec:ordeg}

\subsection{Orientation}
As it matters here, we stress that the Tate twist $\hat{\bZ}(1)$ has an anabelian definition as the geometric fundamental group of $\Gm$.
An \textbf{orientation} on $\pi_1(U/k)$ consists of an isomorphism $\tau$ of the orientation module $\orient_{\pi_1(U/k)}$ with $\hat{\bZ}(1)$. From the anabelian local cohomology sequence (\ref{eq:anabloc}) and the comparison with (\ref{eq:loc}) we see that an orientation on $\pi_1(U/k)$ induces canonically \textbf{local orientations} at each cusp $\ty$, which by definition is a family of  isomorphisms
$\tau_\ty : \rI_\ty \to \hat{\bZ}(1)$ such that the following holds.
\begin{itemize}
\item[(i)] Equivariance with respect to $\pi_1(U)$ and the cyclotomic character $\chi^\cyl$ of $\Gal_k$: for $\gamma \in \pi_1(U)$ and a cusp $\ty$ the following commutes
\[
\xymatrix{  \rI_\ty \ar[d]^{\tau_\ty} \ar[r]^{\gamma()\gamma^{-1}} & \rI_{\gamma.\tilde{y}/\gamma.y} \ar[d]^{\tau_{\gamma.\tilde{y}/\gamma.y}} \\
\hat{\bZ}(1)  \ar[r]^{\chi^\cyl(\gamma) \cdot} & \hat{\bZ}(1).
}
\]
\item[(ii)] The kernel of the map $\bigoplus_{y \in \cp(U)}  \rI_y \to \pig_1^\ab(U)$
from (\ref{eq:anabloc}) becomes the diagonally embeded $\hat{\bZ}(1)$ under the identification coming from the local orientations, which is well defined by (i).
\end{itemize}
We can view local orientations as a $\pi_1(U)$ equivariant map $\tau : W_{-2}(U) \to \hat{\bZ}(1)$ that satisfies a global consistency condition.

\subsection{Standard orientation} The standard orientation on $\pi_1(U/k)$ is given locally by the tame characters and globally by the evaluation at the fundamental class
\[ \orient_{\pi_1(U/k)} = \rH^2(X \otimes k^\alg,\bQ/\bZ)^\vee = \Hom\big(\rH^2(X \otimes k^\alg,\hat{\bZ}(1)), \hat{\bZ}(1)\big) \xrightarrow{\cong} \hat{\bZ}(1),\]
whichever is applicable for $\pi_1(U/k)$. The core of the argument in section \ref{subsec:comp} showed that tame character and evaluation at the fundamental class are indeed compatible in cases both definitions make sense. 

An arbitrary orientation on $\pi_1(U/k)$ will be a multiple of the standard orientation by a factor $\ep \in \hat{\bZ}^\ast$. We call the resulting orientation `\textit{$\ep$ times the standard orientation}' or simply the $\ep$ orientation.

\subsection{Degree} The definition of a degree depends on a choice of orientation. We will silently choose the standard orientation in the sequel.

Let $\ph:\pi_1(U/k) \to \pi_1(V/k)$ be a weight preserving map of fundamental groups of hyperbolic curves $U$ (resp.\ $V$) over $k$ which are complements of divisors $Y$  (resp.\ $Z$) on geometrically connected, smooth projective curves $X$  (resp.\ $W$).  Taking into account orientations the induced map 
$\orient_{\pi_1(U/k)} \to \orient_{\pi_1(V/k)} $
becomes a map $\hat{\bZ}(1) \to \hat{\bZ}(1)$ which is multiplication  by $\deg{\ph} \in \hat{\bZ}$, the \textbf{degree}  of $\ph$.
%Question: How to forbid ell part of deg(ph) = 0? 

Let $\tilde{y} \in \tcp(U)$ be a cusp with $\ph(\rI_\ty) \not= 1$ and $\tilde{z} \in \tcp(V)$ the unique cusp such that $\ph(\rI_\ty) \subset \rI_{\tilde{z}/z}$. The \textbf{local degree} or 
\textbf{ramification index} of $\ph$ at $\tilde{y}$ over $\tilde{z}$ is the unique $\ep_\ph(\tilde{y}/\tilde{z}) \in \hat{\bZ}$, so that $\ph|_{\rI_\ty} : \rI_\ty \to \rI_{\tilde{z}/z}$ equals multiplication by $\ep_\ph(\tilde{y}/\tilde{z})$ after identifying the inertia groups with $\hat{\bZ}(1)$ according to the chosen local orientations. As the restriction to $W_{-2}(U)$ of $\ph$ is equivariant with respect to conjugation, the local degree only depends on the cusps $y \in \cp(U)$ and $\ph(y) := z \in \cp(V)$ and we write $\ep_\ph(y/z) = \ep_\ph(\tilde{y}/\tilde{z})$.

\subsection{Fundamental equation} The following absorbs an argument exploitet by Pop in the context of birational anabelian geometry.
\begin{prop}
Let $\ph:\pi_1(U/k) \to \pi_1(V/k)$ be a weight preserving map between oriented fundamental groups as above. Then for any cusp $z \in \cp(V)$ the fundamental equation
\[ \deg(\ph) = \sum_{y \mapsto z} \ep_\ph(y/z)\]
holds in $\hat{\bZ}$, where the sum extends over cusps $y \in \cp(U)$ so that $\ph(y) = z$.
\end{prop}
\begin{pro}
The map $\ph$ induces a map between the anabelian local cohomology sequences for $\pi_1(U/k)$ and $\pi_1(V/k)$ which using the orientation yield the following commutative diagram
\[
\xymatrix{ \hat{\bZ}(1) \ar[r]^(0.3)\Delta \ar[d]^{\cdot \deg(\ph)} & \bigoplus_{y \in \cp(U)} \hat{\bZ}(1) \ar[d]^{\cdot \ep_\ph(y/z)} \\
\hat{\bZ}(1) \ar[r]^(0.3)\Delta & \bigoplus_{z \in \cp(V)} \hat{\bZ}(1).}
\]
As the maps $\Delta$ are the diagonal embeddings, the assertion follows immediately.
\end{pro}

If we restrict the discussion to isomorphisms $\ph:\pi_1(U/k) \to \pi_1(V/k)$ preserving the weight filtration, we obtain induced bijections $\ph: \cp(U) \to \cp(V)$ and the fundamental equation shows that $\deg(\ph) = \ep_\ph(y/z) \in \hat{\bZ}^\ast$ for all cusps $y$ and $z=\ph(y)$. 
An \textbf{orientation preserving} isomorphism is such an isomorphism $\ph$ with $\deg(\ph)=1$.

%%%%%%%%%%%%%%%%%%%%%%%%%%%%%%%%%%%%%%%%%%%

\section{Anabelian theory of units}  \label{sec:units}

\subsection{Automorphisms of $\pi_1(\Gm/k)$}
\begin{prop}
The automorphism group of the extension $\pi_1(\Gm/k)$ is isomorphic to the semidirect product 
$\widehat{k^\ast} \rtimes \hat{\bZ}^\ast$ with respect to the natural action of $\hat{\bZ}^\ast$ on pro-$\bN$ completions. The group of automorphisms which act as the identity on $\hat{\bZ}(1)$ is $\widehat{k^\ast}$.
\end{prop}
\begin{pro}
The group $\pi_1(\Gm)$ is a semi direct product with respect to the section $s_1$ induced by the rational point $1 \in \Gm(k)$. We write elements of $\pi_1(\Gm)$ as $a\sigma=as_1(\sigma)$ with $a \in \hat{\bZ}(1)$ and $\sigma \in \Gal_k$. 

A unit $\ep \in \hat{\bZ}^\ast$ gives an automorphism $\ep: a\sigma \mapsto a^\ep \sigma$.
 An element $\alpha \in \widehat{k^\ast}$ corresponds by Kummer theory  to the cocycle on $\Gal_k$ with values in $\hat{\bZ}(1)$ defined as 
 \[ \chi_\alpha : \sigma \mapsto \chi_\alpha(\sigma) =\big(\sigma(\sqrt[n]{\alpha})/\sqrt[n]{\alpha}\big)_{n\in \bN} \in \hat{\bZ}(1).\]
  The automorphism \textbf{`scaling with $\alpha$' } is defined as $\alpha: a\sigma \mapsto a\chi_\alpha(\sigma)\sigma$.
We have $\ep \alpha \ep^{-1} : a\sigma \mapsto \big(a^{\ep^{-1}}\chi_\alpha(\sigma)\big)^\ep \sigma = a\chi_\alpha(\sigma)^\ep \sigma$ which equals the automorphism `scaling by $\alpha^\ep$'. The remaining assertions follow immediately.
\end{pro}

Although one is tempted to think of  the projection to $\hat{\bZ}^\ast$ as associating to an automorphism its degree, this interpretation fails. We cannot decide wether an automorphism fixes or interchanges the two cusps of $\Gm$ which are indistinguishable within $\pi_1(\Gm/k)$.

\begin{prop}
Let $x$ be a parameter for $\Gm$. For $\ep = \pm 1$ and $\alpha \in k^\ast$ the geometric automorphism $x \mapsto \alpha x^\ep$ of $\Gm$ induces the automorphism $(\alpha,\ep) \in \widehat{k^\ast} \rtimes \hat{\bZ}^\ast$ of $\pi_1(\Gm/k)$.
\end{prop}
\begin{pro}
The map $x \mapsto x^\ep$ preserves the section $s_1$ and the effect on $\hat{\bZ}(1)$ can be computed from cohomology, thus multiplies with $\ep$. We find again the automorphism denoted $\ep$ in the proof above.

The scaling map $x \mapsto \alpha x$ acts trivial on \'etale cohomology of $\Gm$ by homotopy invariance, hence the associated map on $\pi_1(\Gm/k)$ has $\ep =1$. It remains to observe how the section $s_1$ is shifted. But this is done using the Kummer character $\chi_\alpha(-)$ which mediates  an isomorphism of fibre functors from fibres above $1$ to fibres above $\alpha$, hence `scaling by $\alpha$' acts on sections as the geometric map $x \mapsto \alpha x$ does.
\end{pro}

\subsection{Units} The anabelian way to define \textbf{units} on $\pi_1(U/k)$ is by morphisms to 
$\pi_1(\Gm/k)$. The extension $\pi_1(\Gm/k)$ is a group object with multiplication
\[ \pi_1({\rm mult}) : \pi_1(\Gm) \times_{\Gal_k} \pi_1(\Gm) = \pi_1(\Gm \times_k \Gm) \to \pi_1(\Gm)\]
and unit $s_1 : \Gal_k \to \pi_1(\Gm)$. Hence we obtain a group 
\[ \hat{\OO}^\ast\big(\pi_1(U/k)\big) := \Hom\big(\pi_1(U/k),\pi_1(\Gm/k)\big),\]
which we call the group of \textbf{pro-units} on $\pi_1(U/k)$.  The \textbf{constant pro-units} are those which factor through $\Gal_k$.

If $\pi_1(U/k)$ is equipped with an orientation, then the group $\OO^\ast\big(\pi_1(U/k)\big)$ of \textbf{units} on $\pi_1(U/k)$ is defined as the subgroup of $\hat{\OO}^\ast\big(\pi_1(U/k)\big) $ of maps  that for each cusp $y \in \cp(U)$ induce integral maps $\rI_y = \hat{\bZ}(1) \to \hat{\bZ}(1)$, i.e., the  multiplication by some $\ep_y \in \bZ \subset \hat{\bZ}$.
\begin{prop} 
There is a natural isomorphism $\hat{\OO}^\ast\big(\pi_1(U/k)\big) = \rH^1\big(\pi_1(U),\hat{\bZ}(1)\big)$, under which 
\begin{itemize}
\item[(1)]  constant pro-units  form the image of inflation $ \rH^1\big(k,\hat{\bZ}(1)\big) \inj \rH^1\big(\pi_1(U),\hat{\bZ}(1)\big)$, and
\item[(2)] the natural map $\kappa: \OO^\ast(U) \to \hat{\OO}^\ast\big(\pi_1(U/k)\big)$ given by $f \mapsto \kappa_f := \pi_1(f)$ corresponds to the boundary operator $\delta_\kum$ of Kummer theory on $U$.
\end{itemize}
\end{prop}
\begin{pro}
A unit is nothing but an equivalence class of sections of the extension $\pi_1(\Gm/k)$ pulled back via 
$\pi_1(U) \to \Gal_k$. The part (1) proposition follows from Appendix \ref{app:h1}.

Let $\chi_f: \pi_1(U) \to \hat{\bZ}(1)$ be the Kummer character $\gamma \mapsto \chi_f(\gamma)= (\gamma(\sqrt[n]{f})/\sqrt[n]{f})_{n \in \bN}$, and let $\pr:\pi_1(U) \to \Gal_k$ be the projection.
For (2) and $\gamma \in \pi_1(U)$ we compute 
\[ \kappa_f(\gamma) = \gamma|_{k^\alg(f^\bQ)} = \chi_f(\gamma) \cdot \pr(\gamma) \in \hat{\bZ}(1) \rtimes \Gal_k = \pi_1(\Gm), \]
 which under the identification of (1) becomes $\delta_\kum(f) = \chi_f$.  
\end{pro}

Let $W_{-2}\big(\pig_1^\ab(U)\big) \subset \pig_1^\ab(U)$ be the subgroup generated by inertia. The group $ \Div_Y^0(X)$ of divisors of degree $0$ on $X$ with support in $Y$ has a map $\tau$ to $\Hom_k\big(W_{-2}(\pig_1^\ab(U)),\hat{\bZ}(1)\big) $ which assigns to a divisor the corresponding linear combination of tame characters, and which induces an isomorphism 
\[   \Div_Y^0(X) \otimes_\bZ \hat{\bZ} \xrightarrow{\cong} \Hom_k\big(W_{-2}(\pig_1^\ab(U)),\hat{\bZ}(1)\big).
\]
We define a homomorphism $\gamma: \Pic_X^0(k)  \to \Ext^1_k\big(\pig_1^\ab(X),\hat{\bZ}(1)\big)$
 as the composite of the boundary operator of Kummer theory for $\Pic_X^0$, the map $\pig_1(\Pic^0_X) \to \cHom(\pig_1(\Alb_X),\hat{\bZ}(1))$ given by the Weil pairing, and the edge map in the $\Ext$-spectral sequence:
\[ \Pic_X^0(k) \xrightarrow{\delta_\kum} \rH^1\big(k,\pig_1(\Pic^0_X)\big) \xrightarrow[{\rm Weil}] {\cong}
\rH^1\big(k,\cHom(\pig_1(\Alb_X),\hat{\bZ}(1))\big) \xrightarrow[{\rm edge}]{\cong}  \Ext^1_k\big(\pig_1^\ab(X),\hat{\bZ}(1)\big).
\]
\begin{prop} \label{prop:units}
The maps $\kappa, \tau$ and $\gamma$ fit into the following map of exact sequences  
\[
\xymatrix@C-1ex{ 1 \ar[r] \ar[d] & \OO^\ast(U)/k^\ast  \ar[r]^{\divisor} \ar[d] ^{\kappa} & \Div_Y^0(X)  \ar[r] \ar[d]^{\tau} & \Pic_X^0(k)  \ar[d]^{\gamma} \\
\big(T\Pic^0_X\big)(k) \ar@{^(->}[r] &
 \hat{\OO}^\ast\big(\pi_1(U/k)\big) / \widehat{k^\ast} \ar[r]^(0.4){\res}  & 
 \Hom_k\big(W_{-2}(\pig_1^\ab(U)),\hat{\bZ}(1)\big) 
 \ar[r]^(0.55)\delta & \Ext^1_k\big(\pig_1^\ab(X),\hat{\bZ}(1)\big),
}
\]
\end{prop}
\begin{pro}
The bottom row is exact as the long exact $\Ext$-sequence of $\Gal_k$ modules for 
\begin{equation} \label{eq:w}
 0 \to W_{-2}\big(\pig_1^\ab(U)\big) \to \pig_1^\ab(U) \to \pig_1^\ab(X) \to 0,
\end{equation}
using the Weil-pairing for $\Hom_k\big(\pig_1^\ab(X),\hat{\bZ}(1)\big) = \big(T\Pic^0_X\big)(k) $ and the exact sequence of low degree terms
\[  1 \to \hat{\OO}^\ast\big(\pi_1(U/k)\big) / \widehat{k^\ast}  \xrightarrow{\cong}  \rH^1\big(\pig_1(U),\hat{\bZ}(1)\big)^{\Gal_k} \xrightarrow{0} \rH^2\big(k,\hat{\bZ}(1)\big) \xrightarrow{d_2^{0,1}} \rH^2\big(\pi_1(U),\hat{\bZ}(1)\big) 
\]
 from the Hochschild--Serre spectral sequence.  Indeed, the differential $d_2^{0,1}$ is injective, as the projection $\pi_1(U) \to \Gal_k$ is split over an open subgroup and $\rH^2(k,\hat{\bZ}(1)) = \Hom(\bQ/\bZ, \Br(k))$ contains no torsion.
 
 The middle facet is commutative as for a function $f \in \OO^\ast(U)$ the restriction $\kappa_f|_{\rI_y}$ to inertia at $y$ equals $\nu_y(f)$ times the tame character, where $\nu_y$ is the valuation associated to $y \in X$.
 
It remains to show commutativity of the right facet. 
 Let $D \in \Div_Y^0(X)$ be a divisor of degree $0$ supported on $Y$. The extension $E_{D} = \delta(\tau(D))$ is  obtained by pushing the extension of $\Gal_k$ modules (\ref{eq:w}) with $\tau(D)$.
 
 The line bundle associated to $D$ comes from a unique degree $0$ line bundle  on the Albanese variety $A=\Alb_X$ of $X$. The complement $\bL^0(D)$ of the $0$-section in the corresponding geometric line bundle $\bL(D)$ has  an abelian geometric fundamental group $\pig_1(\bL^0(D))$ as its first Chern class vanishes. Hence, from the canonical  trivialisation of $\bL(D)$ above $U$ we obtain a map of extensions of $\Gal_k$ modules
 \[
 \xymatrix@R-1ex{0 \ar[r] & W_{-2}\big(\pig_1^\ab(U)\big) \ar[r]\ar[d]^{\tau(D)} & \pig_1^\ab(U)  \ar[r]\ar[d] & \pig_1^\ab(X) \ar[r] \ar[d]^\cong & 0 \\
 0 \ar[r] & \hat{\bZ}(1) \ar[r] & \pig_1(\bL^0(D)) \ar[r] & \pig_1(A) \ar[r] &  0,}
\]
from which we conclude that  $E_D$ equals $\pig_1(\bL^0(D))$ as an extension of $\Gal_k$ modules. 

Let $i_D : A \to A \times B$ be the inclusion induced by $D$ as a $k$-rational point of $B=\Pic_X^0$. The line bundle $\bL(D)$ is the restriction of the geometric  Poincar\'e bundle  $\bP_{A \times B}$ along $i_D$. We thus find the extension $E_D$ as the image  $E$ of $\pi_1(i_D): \pig_1(\bL^0(D) \inj \pi_1(\bP_{A \times B}^0)$ and can do the corresponding cocycle calculations in $\pi_1(\bP_{A \times B}^0)$. 

The group $\pi_1(\bP_{A \times B}^0)$ is an extension of $\Gal_k$ by $\pig_1(\bP_{A \times B}^0)$ which is a central extension of $\rT A  \times \rT B$ by $\hat{\bZ}(1)$. The factors $\rT A$ and $\rT B$ are Lagrangian subspaces for the associated commutator pairing which  therefore pairs $\rT A$ with $\rT B$ and as such is nothing but the Weil-pairing induced by the Poincar\'e bundle $\bP_{A \times B}$.

Note that the image $E$ of $\pig_1(\bL^0(D))$ does not depend on $D$, whereas the $\Gal_k$ module structure does depend on $D$ being induced by the conjugation action from the image of the full \'etale fundamental group $\pi_1(\bL^0(D))$, or, what amounts to the same, as the conjugation action via 
$s_{0,D} = \pi_1(i_D)\circ s_0: \Gal_k \to \pi_1(\bP_{A \times B}^0)$. Here $s_0$ belongs to $0 \in A(k)$ and $s_{0,D}$ to $(0,D) \in A \times B(k)$.

For $D=0$ the extension $\pig_1(\bL^0(0))$ splits as an extension of $\Gal_k$ modules via $s_{0,0}$. Let $f: \rT A \to E$ be a splitting which is $s_{0,0}(\Gal_k)$ equivariant. 
By Appendix \ref{app:injsec}  Proposition \ref{prop:abkum} the difference cocycle $\sigma \mapsto \Delta_\sigma := s_{0,D}(\sigma)s_{0,0}(\sigma)^{-1} = s_D(\sigma)s_0(\sigma)^{-1}$ represents $\delta_\kum(D)$.  The class 
\[df =  \sigma \mapsto df_\sigma = \big( x \mapsto  \sigma f(\sigma^{-1} ( x))) - f(x)\big)\]
 of $E_D$ in 
$\rH^1\big(k,\cHom(\rT A,\hat{\bZ}(1))\big)$ computes therefore as
\[ df_\sigma (x)  = s_{0,D}(\sigma) f(s_{0,D}(\sigma)^{-1}x \  s_{0,D}(\sigma)) s_{0,D}(\sigma)^{-1}f(x)^{-1} = \Delta_\sigma f(x) \Delta_\sigma^{-1} f(x)^{-1}\]
where the $\Delta_\sigma$ (resp.\ the $f(x)$) lift elements from the second (resp.\ first) component of $\rT A \times \rT B$ to $\pig_1((\bP_{A \times B}^0)$. Consequently, the cocycles $\Delta_\sigma$ and $df_\sigma$ agree under the Weil pairing finishing the proof.
\end{pro}
\begin{cor} \label{cor:g0units}
Assume that $X$ has genus $0$ or $k$ is a finitely generated extension of $\bQ$, and let the orientation fixed on $\pi_1(U/k)$ be $\ep$ times the standard orientation for some $\ep \in \hat{\bZ}^\ast$. Then we have the following equality in $\hat{\OO}^\ast(\pi_1(U/k))$
\[  \OO^\ast(\pi_1(U/k) = \kappa(\OO^\ast(U))^\ep \cdot \widehat{k^\ast}.\]
In particular all units on $\pi_1(U/k)$ correspond to geometric units of $U$ up to an automorphism of $\pi_1(\Gm/k)$.
\end{cor}
\begin{pro} 
The genus $0$ case follows immediately from Proposition \ref{prop:units}. If $k$ is finitely generated over $\bQ$, the Mordell-Weil theorem shows that $\gamma$ is injective and $(T\Pic^0_X)(k)=0$ and the result follows again from Proposition \ref{prop:units}.
\end{pro}

\subsection{Values} A unit $f : \pi_1(U/k) \to \pi_1(\Gm/k)$ can be evaluated in a section $s \in S_{\pi_1(U/k)}$ as follows. The composition $f \circ s$ yields a section of  $\pi_1(\Gm/k)$, hence an element $f(s) \in \widehat{k^\ast}$ which we call the value of $f$ in $s$. Composing $f$ with an automorphism 
$\alpha\ep$ of $\pi_1(\Gm/k)$ changes the value to $\big(\alpha\ep.f\big)(s) = \alpha f(s)^\ep$.

%%%%%%%%%%%%%%%%%%%%%%%%%%%%%%%%%%%%%%%%%%%

\section{Anabelian geometry of genus $0$ curves - revisited}  \label{sec:genus0again}

\subsection{Cuspidal ratios} A \textbf{cuspidal ratio over $k$} on the filtered, oriented fundamental group $\pi_1(U/k)$ is the value $\lambda = \lambda_{f;y,y'} = f(s_y)/f(s_{y'})  \in \widehat{k^\ast}$, where $s_y, s_{y'}$ are cuspidal sections at $k$-rational cusps $y,y' \in Y(k)$ and $f$ is a unit such that $\ker(f)$ contains the inertia subgroups $\rI_y,\rI_{y'}$ at $y$ and $y'$.  
A \textbf{cuspidal ratio} of $\pi_1(U/k)$ is a cuspidal ratio over some finite extension $k'/k$ for the base change $\pi_1(U \otimes k'/k')$ with respect to the induced orientation.

Obviously, the cuspidal ratio $\lambda_{f;y,y'}$  neither depends on the cuspidal section chosen nor changes its value when $f$ is composed with a scaling automorphism of $\pi_1(\Gm/k)$.

As a corollary to Corollary \ref{cor:g0units} we immediately obtain the following.
\begin{prop}
Assume that $X$ has genus $0$ or $k$ is a finitely generated extension of $\bQ$.
If the orientation on $\pi_1(U/k)$ is $\ep$ times the standard orientation for some $\ep \in \hat{\bZ}^\ast$ then $\lambda_{f;y,y'}$ is the $\ep$-power of an element in the image of the natural map $k^\ast \to \widehat{k^\ast}$.
\end{prop}

From now on we will mainly be preoccupied with curves $U$ which are complements of a divisor $Y$ in a smooth, projective curve $X$ of genus $0$. The fundamental group $\pi_1(U/k)$ shall be endowed with its anabelian weight filtration and an orientation.

%--------------------------------------------------------------------------------------------------------------

\subsection{Definition}
We say that the fundamental group $\pi_1(U/k)$ has \textbf{type} $(0,n)$ if the inertia groups generate the geometric fundamental group $\pig_1(U)$, all cusps of $U$ are $k$-rational, and we are given an ordering $y: \{1,\ldots,n\} \to \cp(U)$ sending $i$ to $y_i$.

%--------------------------------------------------------------------------------------------------------------

\subsection{The double ratio}  \label{subsec:dr}
Let $\pi_1(U/k)$ be a fundamental group of type $(0,4)$ and $y:\{1,2,3,4\} \to \cp(U)$ the corresponding ordering of the cusps. We equip $\pi_1(U/k)$ with  $\ep$ times the standard orientation.

Let $j$ be a unit on $\pi_1(U/k)$ such that $\res(j) =  \ep \cdot \tau(y_3 - y_4)$  with the notation as in Proposition \ref{prop:units}. The double ratio of $\pi_1(U/k)$ as a fundamental group of type $(0,4)$ is defined as the cuspidal ratio 
\[ \DV(\pi_1(U/k),y) = \DV(y_1,y_2;y_3,y_4) = \lambda_{j;y_1,y_2}\]
which is the $\ep$ power of an element of the image of $k^\ast \to \widehat{k^\ast}$ and independent of the unit $j$ chosen.

Let $U_\lambda$ be the complement in $\bP^1_k$ of $Y= \{\lambda,1,0,\infty\}$ and equip $\pi_1(U_\lambda)$ with an ordering of the cusps accordingly and $\ep$ times the standard orientation. The resulting fundamental group of type $(0,4)$ has double ratio 
\[ \DV(\pi_1(U_\lambda/k),y) = \lambda^\ep = \DV(\lambda,1;0,\infty)^\ep,\]
where the double ratio $\DV(a_1,a_2;a_3,a_4)$ of a $4$-tuple $(a_i)$ in $k$ is defined by the formula 
\[ \DV(a_1,a_2;a_3,a_4) = \bruch{a_1-a_3}{a_1-a_4} : \bruch{a_2-a_3}{a_2-a_4}.\]
This follows from Corollary \ref{cor:g0units} as $j$ may be chosen to be $\ep \cdot \pi_1(t)$ where $t$ is the coordinate on $\bP^1_k$ and then $j(s_\lambda) = \lambda^\ep$ and $j(s_1) = 1$.

If $\ep =1$ and so $\pi_1(U/k)$ is equipped with the standard orientation then the following two double ratios add up to $1$:
\[  \DV(y_1,y_2;y_3,y_4) + \DV(y_1,y_3;y_2,y_4) = \lambda + (1-\lambda) = 1.\]

%--------------------------------------------------------------------------------------------------------------

\subsection{Characterization of the standard orientation}

A \textbf{geometric subquotient} of a fundamental group $\pi_1(U/k)$ is the quotient of an open subgroup by the normal subgroup generated by a collection of inertia subgroups. A geometric subquotient of type $(0,n)$ is a geometric subquotient together with a choice of a structure of a fundamental goup of type $(0,n)$ on that subquotient. A geometric subquotient of $\pi_1(U/k)$ inherits naturally an orientation from an orientation on $\pi_1(U/k)$. The factor $\ep$ of the orientation remains the same.

\begin{thm} \label{thm:orientation}
Let $U$ be the complement of a divisor $Y$ in a smooth, projective curve $X/k$ of genus $0$. Let $\pi_1(U/k)$ be equiped with $\ep$ times the standard orientation. 

Then we have necessarily $\ep \in \{\pm 1\}$, if for all geometric subquotients of type $(0,4)$ with field of constants a finite extension $k_1/k$ the double ratio is in the image of $k_1^\ast \to \widehat{k_1^\ast}$ and either of the following conditions holds.
\begin{itemize}
\item[(A)] For any pair of prime numbers $p \not=q$ there is a valuation $\nu$ on $k$ which is nontrivial on $F = k \cap \bQ^\alg$ and with value group $\Gamma \subset \bQ$ that is neither divisible by $p$ nor $q$.
\item[(B)] There is a discrete valuation $\nu$ on $k$ with value group $\bZ$ such that $U$ does not have good reduction over the valuation ring of $\nu$, in the sense that in any smooth model of $X$ some cusps of $U$ coalesce.
\end{itemize}

Moreover, if $\ep \in \{\pm 1\}$ and one of the following conditions (C)--(E) holds, and for the above geometric subquotients of type $(0,4)$ suitable preimages  of the double ratios $\DV(y_1,y_2;y_3,y_4)$ and $\DV(y_1,y_3;y_2,y_4)$ in $k_1^\ast$ add up to $1$,  then even $\ep =1$ and thus the standard orientation  receives an anabelian characterization.
\begin{itemize}
\item[(C)]  The field $k$ is a number field with $k\bQ^\ab/\bQ^\ab$ finite.
\item[(D)]  The field $k$ is a function field over $k_0$ and $U$ is not defined over an algebraic extension of $k_0$.
\item[(E)] The map $k_1^\ast \to \widehat{k_1^\ast}$ is injective for every finite extension $k_1/k$.
\end{itemize}
\end{thm}

\begin{rmk}  (1)
Condition (A) is met by  the fields of condition (C), in particular by $\bQ^\ab$, cf. Remark \ref{rmk:Qab}. 
Fields which are finitely generated over $\bQ_p$ satisfy both (A) and (E), and 
condition (B) holds for fields $k$ that satisfy (D). 

(2) If $k$ is a function field over $k_0$ which is algebraically closed and $U = U_0 \otimes_{k_0} k$, then $\pi_1(U) = \pig_1(U) \times \Gal_k$ and the Galois action has nothing to say about the orientation chosen.
\end{rmk}

\begin{pro}
We first deal with case (A). By Belyi's theorem, the fundamental group of the curve $U_\lambda = \bP^1 - \{\lambda,1,0,\infty\}$ over $k(\lambda)$ occurs for all $\lambda \in \bQ^\alg - \{0,1\}$  as a geometric quotient of type $(0,4)$ of $\pi_1(U/k)$, at least if we allow some finite extension $k_1/k$. We obtain $\mu \in k_1^\ast$ with $\lambda^\ep = \mu$ in $\widehat{k_1^\ast}$.

The assumption of the suitable valuations on $k$ transfers to the finite extension $k_1$. Choose a valuation $\nu$ suitable for $p$ and $q$ and set $\gamma = \nu(\lambda)$ and $\delta=\nu(\mu)$. The valuation extends to a map $\widehat{\nu} : \widehat{k_1^\ast} \to \widehat{\Gamma} \to \bZ_p \times \bZ_q$. The relation $\lambda^\ep = \mu$ becomes $\ep \cdot \gamma = \delta$ in $\bZ_p \times \bZ_q$, where moreover $\gamma$ and $\delta$ come from the diagonally embedded $\bQ \cap (\bZ_p \times \bZ_q)$. Thus the $p$ and $q$ component of $\ep$ are rational and agree. Hence $\ep \in \hat{\bZ}^\ast$ is in the diagonally embedded $\bQ \cap \hat{\bZ}^\ast = \{\pm 1\}$. 

Under the assumption (B), we find, upon restriction to a finite field extension $k_1/k$ that makes the cusps of $U$ rational, even a quotient of type $(0,4)$ with double ratio not a unit with respect to an extension of the valuation $\nu$  to $k_1$. The same but simpler argument with $\widehat{\nu} : \widehat{k_1^\ast} \to \hat{\bZ}$ as above shows that $\ep \in \{\pm 1\}$.

For the remaining part we assume that $\ep = -1$ and argue by contradiction. For one of the geometric subquotients of type $(0,4)$ considered above, let $\lambda$ be the double ratio in the geometric sense. Then by assumption we find $a,b$ in the group $(k_1^\ast)_{\divisor}$ of divisible elements of $k_1^\ast$  such that 
\begin{equation} \label{eq:-1}
\bruch{a}{\lambda} + \bruch{b}{1-\lambda} = 1.
\end{equation}
Under assumption (E) we have $a=b=1$ and thus $\lambda = \zeta_6^{\pm 1}$ is a primite $6^{th}$ root of unity, which we can avoid by the ubiquity of subquotients of type $(0,4)$ provided by Belyi's theorem.

Condition (D) implies that $a$ and $b$ are constants forcing $\lambda$ to be a constant. But the assumption that $U$ is not defined over the constants enables us to find a geometric subquotient of type $(0,4)$ with nonconstant double ratio, hence a contradiction.

Finally, let us assume condition (C). The equation (\ref{eq:-1}) is equivalent to 
\[ \lambda^2 + (b-a-1)\lambda + a = 0,
\]
where $a$ and $b$ are roots of unity by Lemma \ref{lem:div} below. It follows, that the archimedian absolute value of $\lambda$ is bounded above (actually by $3$) and again the ubiquity of subquotients of type $(0,4)$ provided by Belyi's theorem leads to a contradiction.
\end{pro}

\begin{lem} \label{lem:div}
Let $k$ be a number field such that $k\bQ^\ab/\bQ^\ab$ finite. Then the group $(k^\ast)_{\divisor}$ of divisible elements in $k^\ast$ is contained in its torsion group of roots of unity $\mu_\infty(k)$.
\end{lem}
\begin{cor}
We have $(\bQ^{\ab,\ast})_{\divisor} = \mu_\infty$.
\end{cor}
\begin{pro}
We argue by contradiction and assume $a \in (k^\ast)_{\divisor}$ is not a root of unity. Choose an algebraic number field $F \subset k$ that contains $a$ and such that $F\bQ^\ab = k\bQ^\ab$. Then there is a prime number $p \geq 3$ such that $\sqrt[p]{a} \not\in F$ and consequently the field $F_\infty = \bigcup_{n \geq 1} F(\sqrt[p^n]{a}) \subset k\bQ^\ab = F\bQ^\ab$ is an abelian $\bZ_p$ extension of $F$. By the group theory of the cyclotomic character, the field $F_\infty$ is already contained in $\bigcup_{n \geq 1} F(\mu_{p^n})$.

Let $a=\alpha^{p^n}$ in $F(\mu_{p^{n+m}})$.  We apply \cite{rubin}  Lemma (5.7), which gives the injectivity of 
\[ F^\ast/(F^\ast)^{p^N} \inj F(\mu_{p^N})^\ast/\big(F(\mu_{p^N})^\ast\big)^{p^N}\]
for $N = n+m$. We find $\alpha_0 \in F$ with $a^{p^m} = \alpha^{p^{n+m}} = \alpha_0^{p^{n+m}}$.
So $a\zeta$ is divisible by $p^n$ in $F$ for a suitable root of unity $\zeta \in F$. As $\mu_\infty(F)$ is finite, one $\zeta$ is sufficient for all $n$ and with this choice we find $a\zeta \in \bigcap_{n} (F^\ast)^{p^n} \subset \mu_\infty(F)$, hence $a$ was torsion itself.
\end{pro}

%--------------------------------------------------------------------------------------------------------------

\subsection{Preserving orientation}

\begin{thm} \label{thm:preserveorientation}
Let $U$ (resp.\ $V$) be the complement of a divisor $Y$ (resp.\ $Z$) in a smooth, projective curve $X/k$ (resp.\ $W$) of genus $0$, such that $U$ and $V$ are hyperbolic. 

If  conditions (C) or (D) from Theorem \ref{thm:orientation}  or condition (E) and at least one of conditions
 (A) or (B) from Theorem \ref{thm:orientation} hold, then any isomorphism 
$ \ph: \pi_1(U/k) \to \pi_1(V/k)$ that respects the anabelian weight filtration automatically preserves the standard orientation.
\end{thm}

\begin{rmk} 
Condition (C) guarantees by Lemma \ref{lem:Ribet} and Corollary \ref{cor:anWF} that we can characterize the anabelian weight filtration, which  therefore is then automatically respected by isomorphisms.
\end{rmk}

\begin{pro}
As the isomorphism $\ph$ preserves the anabelian weight filtration it also preserves the notion of a geometric subqotient of type $(0,4)$. Let $\lambda \in k^\ast$ be the double ratio of a geometric subqotient of type $(0,4)$ of $\pi_1(U/k)$ and $\mu \in k^\ast$ the double ratio of its companion for $\pi_1(V/k)$. Then the equation 
$ \lambda^{\deg(\ph)} = \mu$ holds in $\widehat{k^\ast}$. As the same $\lambda$'s as in Theorem \ref{thm:orientation} are available we may argue as in the proof of Theorem \ref{thm:orientation} and conclude that  $\deg(\ph) = 1$.
\end{pro}

%--------------------------------------------------------------------------------------------------------------

\subsection{Anabelian geometry of genus $0$}

\begin{thm} \label{thm:anab0n}
Let $U$ (resp.\ $V$) be the complement of a divisor $Y$ (resp.\ $Z$) in a smooth, projective curve $X/k$ (resp.\ $W$) of genus $0$, such that $U$ and $V$ are hyperbolic. 

If  conditions (C) or (D) from Theorem \ref{thm:orientation}  or condition (E) and at least one of conditions
 (A) or (B) from Theorem \ref{thm:orientation} hold, then $U$ and $V$ are isomorphic as $k$-curves if and only if there is an isomorphism $\pi_1(U/k) \to \pi_1(V/k)$  that respects the anabelian weight filtration.
\end{thm}
\begin{rmk} (1) Theorem \ref{thm:anab0n} applies in particular to curves of genus $0$ over $\bQ^\ab$ or to nonconstant curves of genus $0$  over a function field. The latter should be useful in applications to birational anabelian geometry for function fields over algebraically closed fields of charactersitic $0$.

(2) We repeat that condition (C) guarantees that we can characterize the anabelian weight filtration, which  therefore is then automatically respected by isomorphisms.
\end{rmk}
\begin{pro}
By Galois descent as established by Nakamura in \cite{naka2}  Theorem 6.1, we may replace $k$ by a finite extension and therefore assume that all cusps are $k$-rational.

We first discuss the case of genus $0$ with $4$ cusps where moreover both curves carry compatible structures as a fundamental group of type $(0,4)$, say $U \cong U_\lambda$ and $V \cong U_\mu$ with notations as above in Section \ref{subsec:dr}. We can then compare various double ratios and find $\lambda = \mu \cdot a$ and $1-\lambda = (1- \mu) \cdot b$ with 
$a,b \in (k^\ast)_{\divisor}$. It follows that $\lambda=\mu$ or $a\not= b$ and 
\begin{equation} \label{eq:dr}
\lambda = \bruch{a-ab}{a-b} = \DV(a,\infty;ab,b) \qquad \mu = \bruch{b-1}{b-a} = \DV(b,\infty;1,a),
\end{equation}
and so $\lambda,\mu$ are double ratios of divisible elements of $k^\ast$. Under condition (E) this is absurd. In case (D) holds,  the divisible elements are constants and thus also $\lambda$ and $\mu$ are constants. And also assumption (C) restricts the possible values of $\lambda$ and $\mu$ considerably by Lemma \ref{lem:drrou} below. 

For case (E) the above discussion was sufficient. For the cases (C) and (D) we replace $U$ (resp.\ $V$) by a suitable  finite \'etale cover $U' \to U$ of genus $0$ (resp.\ the corresponding $V' \to V$) with again all cusps rational after a finite extension $k'/k$ such that at least one of the double ratios of quotients of $\pi_1(U'/k')$ is not a double ratio of divisible elements. This is possible by Lemma \ref{lem:drrou} and Lemma \ref{lem:div} in case (C) and Belyi's theorem, or  by assumption in case (D). 

Then it is sufficent to treat the case of curves of type $(0,5)$, where moreover one double ratio, say $\alpha$ is not a double ratio of divisible elements and therefore has to agree on both sides by the above argument. We may therefore assume that $U=\bP^1-\{0,1,\alpha,\lambda,\infty\}$ and 
$V=\bP^1-\{0,1,\alpha,\mu,\infty\}$. In order to finish the proof we need to argue that $\mu = \lambda$. Let us assume the contrary, then $\mu$ and $\lambda$ are double ratios of divisible elements as in equation (\ref{eq:dr}). We now compare the following double ratios in $\widehat{k^\ast}$:
\[
\bruch{\lambda-\alpha}{-\alpha} = \DV(\lambda,0;\alpha,\infty) =  \DV(\mu,0;\alpha,\infty) = \bruch{\mu-\alpha}{-\alpha}.
\]
So there is a divisible element $c$ such that $\lambda-\alpha = c\cdot(\mu-\alpha)$. Using (\ref{eq:dr}) we obtain $\lambda = \mu$ or $c \not=1$ and 
\[ \alpha = \bruch{\lambda - c\mu}{1-c} = \bruch{c-a}{c-1} : \bruch{b-a}{b-1} = \DV(c,b;a,1), \]
which finally gives a contradiction because $a,b,c,1$ are in  $(k^\ast)_{\divisor}$.
\end{pro}

\begin{lem} \label{lem:drrou}
Let $\lambda = \DV(\zeta_1,\zeta_2;\zeta_3,\zeta_4)$ be a double ratio of roots of unity $\zeta_i$. Then for any $p$ the $p$-adic absolute value of $\lambda$ satisfies $|\lambda|_p \leq 4$.
If moreover $\lambda \in \bQ$, then $\lambda$ is among the following values: 
\[2,-1,1/2; \quad  3, 1/3, -2, -1/2, 2/3,3/2; \quad 4,1/4,-3,-1/3,4/3,3/4.\]
\end{lem}
\begin{pro}
For $\zeta$ a primitive $n^{th}$-root of unity, the value $\zeta-1$ is a unit except if $n=p^m$ and then it is a uniformizer in $\bQ_p(\mu_{p^m})$. We conclude that $1 \geq |\zeta -1|_p \geq p^{-1/p^{m-1}(p-1)} \geq 1/2$, hence
\[   |\lambda|_p =  \left|\DV(\zeta_1,\zeta_2;\zeta_3,\zeta_4)\right|_p = \left|\bruch{\zeta_1-\zeta_3}{\zeta_1-\zeta_4} : \bruch{\zeta_2-\zeta_3}{\zeta_2-\zeta_4}\right|_p  = \left| \bruch{1-\zeta_3\zeta_1^{-1}}{1-\zeta_4\zeta_1^{-1}} : \bruch{1-\zeta_3\zeta_2^{-1}}{1-\zeta_4\zeta_2^{-1}}\right|_p \leq 4  \]
For $\bQ$-rational double ratios of roots of unity we observe that we have a natural $S_3$ action on the set of values including $\lambda \mapsto 1/\lambda$, and that by the estimate above the only prime factors that may occur are $2$ and $3$ with only $2$ possibly occuring twice. Imposing this on the full $S_3$ orbit already pins down  $\lambda$ into the given list, whose values by the way are all attained for suitable roots of unities.
\end{pro}

%%%%%%%%%%%%%%%%%%%%%%%%%%%%%%%%%%%%%%%%%%%

\section{Anabelian geometry of  $\cM_{0,4}$ and of the $j$-invariant}  \label{sec:moduli}

\subsection{}  \label{subsec:m04}
The moduli space $\cM_{0,4}$ parameterizes pairs $(X/S; \underline{y}) \in \cM_{0,4}(S)$ where $X/S$ is a smooth, projective $S$-curve of genus $0$ with $4$ marked points $\underline{y} = y_1,y_2,y_3,y_4 \in X(S)$ with pairwise disjoint images. The double ratio $\DV(\underline{y})$ of a marked curve $(X/S, \underline{y})$ is 
\[\DV(\underline{y}) =  \bruch{y_1-y_3}{y_1-y_4} : \bruch{y_2-y_3}{y_2-y_4} \in \OO^\ast(S),\]
where the values $y_i$ are taken with respect to and are independent of any choice of parameter on $X \cong \bP^1$ locally on $S$. The double ratio map is an isomorphism $\DV: \cM_{0,4} \to \bP^1-\{0,1,\infty\}$.

\subsection{} The map ${\rm SC}$ of the section conjecture for $\cM_{0,4}$ factors as follows.
\[
\xymatrix@M+1ex@C-2ex@R-3ex{ \cM_{0,4}(k) \ar[dd]^{\DV}_{\cong} \ar[rr]^{\rm SC} \ar[dr]_(0.3){\pi_1} && S_{\pi_1(\cM_{0,4}/k)}  \ar[dd]^{\DV}_{\cong} \\
&
\left\{  {\begin{array}{c}
\text{fundamental groups } \pi_1(U/k)
\text{ of type } (0,4) \\
\text{ with  anabelian weight structure}
\end{array}}
 \right\} \ar[rd]^(0.65){\Theta} & \\
 \big(\bP^1 \setminus \{0,1,\infty\}\big)(k)  \ar[rr]^{\rm SC}  && S_{\pi_1(\bP^1_k-\{0,1,\infty\}/k)}
}
\]
The map $\pi_1: (X/k,\underline{y})  \mapsto \pi_1(U/k)$ with $U=X-\{\underline{y}\}$  and equipped with the specified extra structure is surjective by definition and bijective if curves of type $(0,4)$ over $k$ are anabelian in a weak sense.   

The map $\Theta$ has  an anabelian definition as follows. The quotient of $\pi_1(U/k)$ by the normal subgroup generated by the inertia subgroups above the cusp $y_1$ is geometric and thus canonically isomorphic to $\pi_1(\bP_k^1-\{0,1,\infty\}/k)$.  A cuspidal section associated to $y_1$ of $\pi_1(U/k)$ maps to a unique geometric section $\Theta(\pi_1(U/k))$  of   $\pi_1(\bP^1_k-\{0,1,\infty\}/k)$. The diagram commutes, because the curve $\big(\bP^1_k,\{\lambda,1,0,\infty\}\big)$ maps to the section $s_\lambda$ of $\pi_1(\bP^1_k-\{0,1,\infty\}/k)$, where $\lambda=\DV(\lambda,1;0,\infty)$ equals the double ratio of the $k$-point  $\big(\bP^1_k,\{\lambda,1,0,\infty\}\big) \in \cM_{0,4}(k)$.

\subsection{}  The universal curve of type $(0,4)$ is the map $\cM_{0,5} \to \cM_{0,4}$ that forgets the fifth point. We abreviate the fundamental group of a chosen geometric fibre by $\pig_{0,4}$ and obtain an extension
\begin{equation} \label{eq:univ}
1 \to \pig_{0,4} \to \pi_1(U_{\cM_{0,5}}) \to \pi_1(\cM_{0,4}) \to 1,
\end{equation}
where $U_{\cM_{0,5}}$ is the locus of triviality of the regular log structure on $\cM_{0,5}$ defined by the $4$ universal sections.  The pullback along section $y_i$ defines a log structure $M_i$ on $\cM_{0,4}$ with constant rank $1$.  We obtain a diagram of (logarithmic) fundamental groups
\[
\xymatrix{ 1 \ar[r] & \hat{\bZ}(1) \ar[d]^{y_i} \ar[r] & \pi_1^{\log}(\cM_{0,4}, M_i) \ar[d]^{\pi_1(y_i)} \ar[r] & \pi_1(\cM_{0,4}) \ar@{=}[d] \ar[r] & 1 \\
1 \ar[r] & \pig_{0,4} \ar[r] & \pi_1(U_{\cM_{0,5}}) \ar[r] & \pi_1(\cM_{0,4}) \ar[r] & 1,} 
\]
where $y_i$ maps $\hat{\bZ}(1)$ injectively onto an inertia group above $y_i$ and the extension in the top row is the first Chern class extension of $y_i^\ast\big(\OO(y_i)\big)$. The top row splits as $\Pic(\cM_{0,4})=1$. Any splitting leads to a section of the bottom row which cyclotomically normalises the inertia at $y_i$ and thus is a cuspidal section at $y_i$. 

As a converse map to $\Theta$ we define for a section $s \in  S_{\pi_1(\cM_{0,4}/k)}$ the extension $\prod_s$ as the pullback via $s$ of the universal extension (\ref{eq:univ}), which is an  extension of 
$\Gal_k$ by $\pig_{0,4}$ endowed with an anabelian weight filtration and $4$ ordered cuspidal sections. 
It is not difficult to see that  $s \mapsto \prod_s$ restricted to geometric sections and  $\Theta$ are mutually inverse. It is tempting to believe that a more refined study of (\ref{eq:univ}) may even give a proof of the section conjecture for $\cM_{0,4}$.

%-------------------------------------------------------------------------------------------------------------------

\subsection{Proper units} In order to discuss the case of higher genus, in particular elliptic curves, we need to address again Proposition \ref{prop:units}. A cuspidal ratio depends only on the restriction of the unit to its divisor in $\Div^0_Y(X) \otimes_\bZ \hat{\bZ} =  \Hom_k\big(W_{-2}(\pig_1^\ab(U)),\hat{\bZ}(1)\big) $ up to a cuspidal ratio that comes from a unit in the image of 
$\big(T\Pic^0_X\big)(k) \subset \hat{\OO}^\ast\big(\pi_1(U/k)\big) / \widehat{k^\ast} $, which we choose to call \textbf{proper units} and which are constructed as follows. To an $L \in \big(T\Pic^0_X\big)(k)$ belongs by the Weil-pairing a $\Gal_k$ invariant map $\pig_1^\ab(X)  \to \pig_1 \Alb_X \to \hat{\bZ}(1)$, which defines a map  of semidirect products $\pi_1(\Alb_X/k) \to \pi_1(\Gm/k)$ and thus the proper unit 
\[ f_L : \pi_1(U/k) \to \pi_1(X/k) \to \pi_1(\Alb_X/k) \to \pi_1(\Gm/k)\]
which has trivial divisor. Using Kummer theory for $\Alb_X$ and Appendix \ref{app:injsec} Proposition \ref{prop:abkum}, we get a pairing 
\[ \big(T\Pic^0_X\big)(k) \times \widehat{\Alb_X(k)} \to \widehat{k^\ast} \] 
that maps $(L,y'-y)$ to the cuspidal ratio $\lambda_{f_L;y',y}$. This immediately yields the following lemma.
\begin{lem} \label{lem:valueproperunits}
The cuspidal ratio of a proper unit with respect to a pair of points which differ by  $n$-torsion in $\Alb_X$ takes values in $n$-torsion of $\widehat{k^\ast}$.
\end{lem}
\begin{lem} \label{lem:torsion}
The $n$-torsion of $\widehat{k^\ast}$ is generated by the image of $\mu_n(k)$.
\end{lem}
\begin{pro}
Let $a \in  \widehat{k^\ast}$ with representatives $a_r$ mod $(k^\ast)^r$ be $n$-torsion, i.e., for $r\in \bN$ we have $\alpha_r \in k^\ast$ such that $a_{r+n}^n = \alpha_r^{r+n}$. For some $\zeta \in \mu_n(k)$ depending on $r$ we thus have $a_{r+n}\zeta = \alpha_r^r$. As there are only finitely many $n^{th}$ roots of unity, one $\zeta$ will do it for all $r$. Then $a_r \zeta \equiv a_{r+n}\zeta \equiv 1$ modulo $(k^\ast)^r$ shows that $a=\zeta^{-1}$.
\end{pro}

\subsection{The $j$-invariant}\footnote{The author acknowledges a discussion with Hiroaki Nakamura held in 2005 where the question how to recover the $j$-invariant was addressed while the author enjoyed the hospitality of the University of Okayama.}
Let $E^0 = E - \{e\}$ be the complement of the origin in an elliptic curve $E/k$. We consider $\pi_1(E^0/k)$ as equipped with the anabelian weight filtration and a choice of cuspidal section $s_e$ at  $e$. The finite \'etale cover $E-E[2] \to E^0$ which is given by multiplication by $2$ is characterized as the unique $V_4$ cover which is unramified over $E$ together with a choice of a lift of  $s_e$. The set of cusps of $E-E[2]$ as a $\Gal_k$ set thus has a distinguished element, again denoted $e$ and three other elements $P_0,P_1$ and $P_\infty$. Let $k'/k$ be the smallest Galois extension such that all cusps of $E-E[2]$ are rational over $k'$.  We get a homomorphism
\[ \Gal(k'/k)  \to S_3 ={\rm Perm}\{0,1,\infty\} , \]
which sends $\sigma \mapsto \pi_\sigma$ with $\sigma(P_i) = P_{\pi_\sigma(i)}$. Moreover, we have a homomorphism 
\[g: S_3 ={\rm Perm}\{0,1,\infty\}  \to \Aut(\bP^1),\]
 which sends a permutation $\alpha$ to the unique automorphism $g_\alpha$ of $\bP^1$ which permutes $0,1, \infty$ accordingly, e.g.,  $g$ maps the transpositions $(0,\infty)$ and $(0,1)$ to the automorphisms $t \mapsto 1/t$ and $t \mapsto 1-t$ respectively.

Let $x_i$ be a unit on $\pi_1(E-E[2]/k)$ with divisor $2P_i - 2e$. For any permutation $\alpha \in S_3$ we define the cuspidal ratio
\[ \lambda_\alpha = \lambda_{x_{\alpha(1)}; P_{\alpha(\infty)}, P_{\alpha(0)}} =  x_{\alpha(1)}\big(P_{\alpha(\infty)}\big)/ x_{\alpha(1)}\big(P_{\alpha(0)}\big) \in k'^\ast/\pm k'^\ast_{\divisor}, \]
which is well defined by Lemma \ref{lem:valueproperunits} and Lemma \ref{lem:torsion}. 
Fixing the points of order $2$ imposes on $E$ the structure of a twist of the Legendre form as 
\[E=\{ cy^2=x(x-1)(x-\lambda)\}\]
 such that  $e$ is the point at infinity and $P_0=(1,0), P_1=(0,0)$ and $P_\infty = (\lambda,0)$ for some $\lambda \in k'^\ast$ different form $1$.  We then have  $x_0= \pi_1(1-x)$, $x_1 = \pi_1(x)$ and $x_\infty = \pi_1(\lambda-x)$ unique up to a multiplicative constant and a proper unit. 
We find that 
\[g_\alpha(\lambda) \equiv \lambda_{\alpha^{-1}}  \text{ modulo }  \pm k^\ast_{\divisor}.\]
Let us assume that $\mu \in k'^\ast -\{1\}$ also gives rise to these values $\lambda_\alpha$. Then there are elements $\zeta, \xi \in \pm k'^\ast_{\divisor}$ such that $\lambda = \zeta \mu$ and $1-\lambda = \xi(1-\mu)$. Therefore either  $\lambda$ equals $\mu$ or $\zeta \not= \xi$ and 
\[ \lambda= \bruch{\zeta - \zeta \xi}{\zeta-\xi} = \DV(\zeta^{-1},\infty;1,\xi^{-1})
\qquad
\mu = \bruch{1-\xi}{\zeta-\xi} = \DV(\xi,\infty;1,\zeta) 
\]
are double ratios of elements from $\pm k'^\ast_{\divisor} \cup \{\infty\}$. Whenever the latter turns out not to be the case, we may define the $j$-invariant of the fundamental group of a curve of type $(1,1)$ by the usual formula as
\[ j(\pi_1(E^0/k)) = j(\lambda) = 1728 \bruch{(\lambda^2 - \lambda + 1)^3}{\lambda^2(1-\lambda)^2} \in k ,\]
which is  independent of the choice of $\lambda$ as $j(\lambda)$ is invariant under the given $S_3$ action and  in fact gives a value in  $k$ as by structure transport 
$\sigma(\lambda_\alpha) = \lambda_{\pi_\sigma \circ \alpha}$ and so the $S_3$ orbit is $\Gal(k'/k)$ invariant.  Consequently, we obtain the following weak anabelian result.  
\begin{thm}
Let $E/k$ be an elliptic curve with origin $e \in E(k)$. Assume that $k^\ast \inj \widehat{k^\ast}$ is injective or that $k$ is a function field over $k_0$ and $E$ is not defined over $k_0$.  

Then the above construction decodes grouptheoretically the $j$-invariant of $E$ from the fundamental group extension $\pi_1(E-\{e\}/k)$.
\end{thm}

%%%%%%%%%%%%%%%%%%%%%%%%%%%%%%%%%%%%%%%%%%%
%%%%%%%%%%%%%%%%%%%%%%%%%%%%%%%%%%%%%%%%%%%

\begin{appendix}

%%%%%%%%%%%%%%%%%%%%%%%%%%%%%%%%%%%%%%%%%%%
%%%%%%%%%%%%%%%%%%%%%%%%%%%%%%%%%%%%%%%%%%%

\section{Torsors, sections and non-abelian $\rH^1$} \label{app:h1}

\subsection{}
Let $\Gamma$ be a pro-finite group. A \textbf{$\Gamma$-group} is a pro-finite group $N$ together with a continuous action of $\Gamma$ by group automorphisms.

\subsection{}
A \textbf{$1$-cocycle} of $\Gamma$ with values in $N$ is a continuous map $a: \Gamma \to N$ such that 
\[ a_{\sigma\tau} = a_\sigma \sigma(a_\tau)\]
for all $\sigma,\tau \in \Gamma$. Two $1$-cocycles $a,b$ are by definition equivalent if  there is a $c \in N$ such that for all $\sigma \in \Gamma$ we have 
\[ b_\sigma = c^{-1} a_\sigma \sigma(c). \]
The \textbf{first non-abelian cohomology} $\rH^1(\Gamma,N) $ of $\Gamma$ with values in $N$ is the pointed set of equivalence classes of continuous $1$-cocycles with the trivial cocyle $a_\sigma \equiv 1$ as the special element.

\subsection{}
The \textbf{semidirect product} $N \rtimes \Gamma$ sits by definition in a canonically split short exact sequence 
\[ 1 \to N \to N \rtimes \Gamma \to \Gamma \to 1 \]
such that for $n \in N$ and $\sigma \in \Gamma$ we have $\sigma n = \sigma(n) \sigma$, when $\sigma$ is considered as an element of $N \rtimes \Gamma$ via the canonical splitting. Let $S_{N \rtimes \Gamma}$ denote the pointed set of $N$-conjugacy classes of sections of $N \rtimes \Gamma \surj \Gamma$ with the canonical splitting as the special element.

\subsection{} 
A \textbf{$\Gamma$-equivariant right $N$-torsor} is a pro-finite set $P$ with a continuous $\Gamma$-action and a continuous free transitive and $\Gamma$-equivariant action by $N$. The pointed set of isomorphy classes of $\Gamma$-equivariant right $N$-torsor will be denote by ${\rm Tors}_\Gamma(N)$ with $P=N$ and right translation as the special element.

\begin{prop}
The pointed sets  $\rH^1(\Gamma,N)$, $S_{N \rtimes \Gamma}$,
and $ {\rm Tors}_\Gamma(N)$ are in natural bijection.
\end{prop}
\begin{pro}
Let $a : \Gamma \to N$ be a continuous map. Then $a$ is a $1$-cocycle if and only if $s_a = (\sigma \mapsto a_\sigma \sigma)$ is a group homomorphism.  If $b_\sigma = c^{-1}a_\sigma \sigma(c)$ for all $\sigma \in \Gamma$ then $s_a = c()c^{-1} \circ s_b$ and conversely. Clearly, every section is of the for $s_a$ for some continuous map $a$ and the trivial cocycle leads to the canonical section.

A $\Gamma$-equivariant right $N$-torsor $P$ is the same as a pro-finite set $P$ with a right $N \rtimes \Gamma$ action such that the subgroup $N$ acts free and transitive. The action is given by $p.(a\sigma) := \sigma^{-1}(pa)$.  The stabilizer $\Gamma_t \subset N\rtimes \Gamma$ of $t \in P$ projects isomorphically onto $\Gamma$ and hence defines a section in $S_{N \rtimes \Gamma}$. Moving the lement $t$ conjugates $\Gamma_t$ by an element of $N$, so that we obtain a well defined map 
$ {\rm Tors}_\Gamma(N) \to S_{N \rtimes \Gamma}$. The inverse maps a section $s$ to the set $P=s(\Gamma)\backslash(N \rtimes \Gamma)$ where the $\Gamma$-equivariant $N$ right torsor structure comes from right translation with $N \rtimes \Gamma$.
\end{pro}

Combining the above constructions we get a bijection $ {\rm Tors}_\Gamma(N) \to \rH^1(\Gamma,N)$ as follows. An element $t \in P$ defines a $1$-cocycle $a^{P,t}:\Gamma \to N$ by the formula $\sigma(t) = t a_\sigma^{P,t}$. Moving $t$ replaces $a^{P,t}$ by an equivalent cocycle. Conversely, a $1$-cocylce $a$ allows to twist the trivial $\Gamma$-equivariant torsor $N$ to a torsor $^aN$ which equals $N$ as a right $N$ torsor but has $\sigma._a n = a_\sigma n$ as twisted $\Gamma$-action.

\subsection{}
Now let us assume that $N$ is an abelian group. This makes $\rH^1(\Gamma,N)$ an abelian group. Let $E$ be an extension of $\Gamma$ by $N$ which splits, the set of $N$-conjugacy classes of sections denoted by $S_E$. An element $s \in S_E$ may now be twisted by a cohomology class $a \in \rH^1(\Gamma,N)$ by choosing a representing cocyle and the formula $a.s(\sigma) := a_\sigma s(\sigma)$. The conjugacy class of $a.s$ is independent of the choice.

A $\Gamma$-equivariant right $N$-torsor $P$ may also be twisted to a torsor $^aP$ by a cohomology class $a \in \rH^1(\Gamma,N)$ by choosing a representing cocyle and the formula for the twisted $\Gamma$ action being $\sigma._a(p) = \sigma(p)a_\sigma$. The isomorphy class of the twist does not depend on the choice.

\begin{prop} \label{prop:h1}
Twisting as above defines free and transitive actions of $\rH^1(\Gamma,N)$ on $S_E$ respectively ${\rm Tors}_\Gamma(N)$.
\end{prop}
\begin{pro}
This is immediate from the above constructions by noting that $s_a = a.s$ and $^aN$ above is in fact really the twist of $N$ by $a$. The assumption on $N$ being abelian is necessary to define the action of $\rH^1(\Gamma,N)$ in the first place.
\end{pro}

%%%%%%%%%%%%%%%%%%%%%%%%%%%%%%%%%%%%%%%%%%%

\section{Injectivity in the section conjecture}  \label{app:injsec}

\subsection{Kummer theory for abelian varieties}
Let $A/k$ be an abelian variety and $n$ invertible in $k$. We consider the Kummer sequence $0 \to A[n] \to A \to A \to 0$ and the truncated mod $n$ fundamental group extension $\pi_1^{[n]}(A/k)$
\[ 1 \to A[n] \to \pi_1^{[n]}A \to \Gal(k^\alg/k) \to 1,\]
which is a quotient of $\pi_1(A/k)$. Hence rational points $a\in A(k)$ lead to sections of $\pi_1(A/k)$ which project onto sections $s_a$ of $\pi_1^{[n]}(A/k)$.
\begin{prop} \label{prop:abkum}
The boundary of the Kummersequence 
\[\delta : A(k)/nA(k) \inj \rH^1(k,A[n])\]
maps $a \in A(k)$ to the class of the difference cocycle $\sigma \mapsto s_a(\sigma) s_0(\sigma)^{-1}$ of the section $s_a$ and the section $s_0$ associated to the origin.
\end{prop}
\begin{pro}
It is sufficient to compute and compare the actions of $\Gal(k^\alg/k)$ on the fibre $A[n](k^\alg)$ of multiplication by $n$ over $0$ that are induced by $s_a$ and $s_0$ respectively. Via $s_0$ the action is the natural one.

The same holds for $s_a$ for the action on the fibre over $a$. In order to regard $s_a$ as a section of $\pi_1(A,0)$, so we can let it act on the fibre above $0$, we have to choose a path from $0$ to $a$. The path is determined by what it does for the cover `multiplication by $n$' and thus is given by the translation by a fixed $n$th root $\bruch{1}{n}a$ of $a$. 

For a point $P \in A[n](k^\alg)$ we compute
\[ s_a(\sigma)\big(P\big) = - \bruch{1}{n}a + \sigma\big(P + \bruch{1}{n}a\big) = - \bruch{1}{n}a + \sigma\big(P\big) + \sigma\big(\bruch{1}{n}a\big) =  \sigma\big(\bruch{1}{n}a\big) - \bruch{1}{n}a + s_0(\sigma)\big(P\big).\]
Hence $s_a(\sigma)s_0(\sigma)^{-1}$ equals translation by $\sigma(\bruch{1}{n}a) - \bruch{1}{n}a$, which is nothing but the value of a cocycle for $\delta(a)$ at $\sigma$.
\end{pro}

\subsection{Subabelian varieties}

\begin{thm} \label{thm:injsec}
Let $k$ be a field such that  the group of $k$-rational poins $A(k)$ has no nontrivial divisible elements for any abelian variety $A/k$.

Let $x,y \in X(k)$ be $k$-rational points of a variety $X/k$ that immerses into an abelian variety such that the associated sections $s_x,s_y \in S_{\pi_1(X/k)}$ agree. Then $x$ equals $y$.  
\end{thm} 

\begin{rmk} 
The assumptions of Theorem \ref{app:injsec} are valid for number fields (resp.\ local fields) by the Mordell-Weil Theorem (resp.\ compactness). In particular the injectivity part of the section conjecture holds for hyperbolic curves over such fields as the Albanese map is eventually an immersion for suitable finite \'etale covers.
\end{rmk}

\begin{pro}
By functoriality of sections we may assume that $X$ is in fact an abelian variety $A/k$. 
For each $n \in \bN$ the sections $s_x, s_y$ induce the same section of $\pi_1^{[n]}(A/k)$. Hence  the difference cocycle $\sigma \mapsto s_x(\sigma)s_y(\sigma)^{-1}$ is trivial in $\rH^1(k,A[n])$ and $x-y$ is divisible by $n$ in $A(k)$ by Proposition \ref{prop:abkum}.   We conclude that $x=y$ from the assumption on divisible elements in $A(k)$.
\end{pro}

%%%%%%%%%%%%%%%%%%%%%%%%%%%%%%%%%%%%%%%%%%%%

%\section{Tangential base points} \label{app:ts}

%depend only on the normal direction!

%
%%%%%%%%%%%%%%%%%%%%%%%%%%%%%%%%%%%%%%%%%%%
%%%%%%%%%%%%%%%%%%%%%%%%%%%%%%%%%%%%%%%%%%%

\end{appendix}

%%%%%%%%%%%%%%%%%%%%%%%%%%%%%%%%%%%%%%%%%%%
%%%%%% Bibliography %%%%%%%%%%%%%%%%%%%%%%%%%%%%%%
%%%%%%%%%%%%%%%%%%%%%%%%%%%%%%%%%%%%%%%%%%%

%\nocite{*}

%\cleardoublepage
%\phantomsection
%\addcontentsline{toc}{chapter}{Literaturverzeichnis}

%\bibliography{../0service/biblio}

%-------------------------------------------------------------------------

%-------------------------------------------------------------------------

\end{document}